\documentclass[3p,10pt]{elsarticle}

\usepackage{amssymb,amsthm}
\usepackage{latexsym,amsmath,subfigure,color,siunitx}

\usepackage{lineno}

\journal{}

\newcommand{\eps}{\varepsilon}
\newcommand{\epsb}{\varepsilon_\mathrm{b}}

\newcommand{\sigmab}{\sigma_\mathrm{b}}
\newcommand{\set}[1]{\left\{#1\right\}}
\newcommand{\abs}[1]{\left|#1\right|}
\newcommand{\e}{\mathrm{e}}
\newcommand{\E}{\mathrm{E}}
\newcommand{\M}{\mathrm{M}}
\newcommand{\N}{\mathrm{N}}

\newcommand{\mE}{\mathbf{E}}

\newcommand{\mH}{\mathbf{H}}
\newcommand{\mU}{\mathbf{U}}
\newcommand{\mV}{\mathbf{V}}
\newcommand{\mW}{\mathbf{W}}
\newcommand{\p}{\partial}
\newcommand{\ma}{\mathbf{a}}

\newcommand{\mf}{\mathbf{f}}

\newcommand{\mr}{\mathbf{r}}
\newcommand{\vt}{\boldsymbol{\theta}}

\DeclareMathOperator*{\inc}{inc}
\DeclareMathOperator*{\tot}{tot}
\DeclareMathOperator*{\area}{area}
\DeclareMathOperator*{\scat}{scat}
\DeclareMathOperator*{\noise}{noise}
\DeclareMathOperator*{\tm}{TM}
\DeclareMathOperator*{\dm}{DM}
\DeclareMathOperator*{\rback}{rb}
\DeclareMathOperator*{\rstar}{r\star}

\theoremstyle{plain}
\newtheorem{thm}{Theorem}[section]

\newtheorem{cor}[thm]{Corollary}

\theoremstyle{remark}
\newtheorem{rem}{Remark}[section]
\newtheorem{ex}{Example}[section]

\begin{document}

\begin{frontmatter}



\title{Application of MUSIC algorithm in real-world microwave imaging of unknown anomalies from scattering matrix\tnoteref{Seo}}
\tnotetext[Seo]{Dedicated to Professor Jin Keun Seo on the Occasion of His 60th Birthday.}

\author{Won-Kwang Park}
\ead{parkwk@kookmin.ac.kr}
\address{Department of Information Security, Cryptography, and Mathematics, Kookmin University, Seoul, 02707, Korea.}

\begin{abstract}
  In this contribution, we consider MUltiple SIgnal Classification (MUSIC)-type algorithm for a non-iterative microwave imaging of small and arbitrary shaped extended anomalies located in a homogeneous media from scattering matrix whose elements are scattering parameters measured at dipole antennas. In order to explain the feasibility of MUSIC in microwave imaging, we investigate mathematical structure of MUSIC by establishing a relationship with an infinite series of Bessel function of integer order and antennas setting. This is based on the representation formula of scattering parameters in the presence of small anomalies and the application of Born approximation. Simulation results using real-data at $f=925$MHz of angular frequency are exhibited to show the feasibility of designed algorithm and to support investigated structure of imaging function.
\end{abstract}


\end{frontmatter}




\section{Introduction}\label{sec:1}
Microwave imaging for identification of the shape or location of unknown anomaly from measured scattered fields or scattering parameters is an interesting and important inverse scattering problem. Throughout various researches, it turned out that this is inherently a non-linear and ill-posed problem so that it cannot be successfully resolved.

In order to solve this interesting and important problem, various inversion algorithms have been investigated for example, Newton-type scheme \cite{SBSSNST} for brain imaging, conjugate gradient \cite{BSSSNST} and Gauss-Newton methods \cite{ML2} for breast imaging, back-projection method \cite{BBPALH} for heart imaging, Newton-Kantorovich method \cite{MJBB} for arm imaging, Levenberg-Marquadt method \cite{CJGT} for in-body imaging, least-mean square approach \cite{TTC} for ground penetrating radar (GPR), linear inversion technique \cite{HSM2} for thermal therapy monitoring, focusing method \cite{KJFF} for damage detection of civil structures, and level-set method \cite{ADIM} for crack reconstruction. We also refer to \cite{CMP,CZBN,FFK,MFLPP,PBCID,SKVH,SKL,SSKLJ} for various application and inversion techniques in microwave imaging.

As we seen, most of algorithms for solving inverse scattering problem related to microwave imaging are based on the Newton-type iteration scheme, i.e., obtaining the shape of the target (minimizer), which minimizes the discrete norm (usually, $\ell^2-$norm) between the measured scattered field or scattering parameter data in the presence of true and man-made targets such that
\[\mbox{find }\mathbf{x}\mbox{ to minimize }\mathcal{E}(\mr):=\frac12||\mathbb{F}\mathbf{x}-\mathbb{M}||^2+\mbox{regularization},\]
where $\mathbb{F}$ is the discretized linear forward operator, $\mathbb{M}$ contains measurements data, and $\mathbf{x}$ is a complex vector of pixels of the differential object. Iteration-based techniques have been applied to identify the shape, location, and topological properties of unknown target. Nevertheless, obtaining a good initial guess, estimating \textit{a priori} information, selecting appropriate regularization term, and evaluating complex Fr{\'e}chet derivative at each iteration steps must be considered beforehand. If one of these conditions is not fulfilled, one may encounter various problems such as the non-convergence, the occurrence of local minimizer, and the requirement of large computational costs due to the large number of iteration procedures. Due to this reason, application of Newton-type method is very hard to the real-world inverse problems.

Correspondingly, for an alternative, various non-iterative inversion techniques have been developed and successfully applied to various inverse scattering problems for example a variational algorithm based on Fourier inversion \cite{AIM,AMV}, direct sampling method \cite{IJZ1,IJZ2}, orthogonality sampling method \cite{P1,CAP}, linear sampling method \cite{CC,KR}, factorization method \cite{K1,K5}, Kirchhoff and subspace migrations \cite{AGKPS,P-SUB3}, topological derivatives \cite{P-TD3,P-TD5}, truncated singular value decomposition (TSVD) inversion scheme \cite{SAS,SSA}. Among them MUltiple SIgnal Classification (MUSIC) algorithms operated at single and multiple frequency have been successfully applied to the various inverse scattering problems for example, detection of inhomogeneities in half-space problem \cite{AIL1,IGLP}, multi-layered medium \cite{SCC}, and inhomogeneous media \cite{K1}, imaging of crack-like defects \cite{AKLP,P-MUSIC1}, detecting internal corrosion \cite{AKKLV}, eddy-current nondestructive evaluation \cite{HLD}, electrical breast biopsy \cite{S2}, imaging of extended targets \cite{HSZ1}, and microwave imaging \cite{PKLS}. We also refer to \cite{RSAAP,RSCGBA} for theoretical and experimental validation of MUSIC in breast cancer detection in monostatic measurement configuration and \cite{AILP,AK2,C,CZ,ML,OBP,P-MUSIC3,PL1,PL3,SL} for applications of MUSIC in various inverse scattering problems.

Following the traditional results, most of non-iterative imaging algorithms were established under some strong assumptions that every elements of so-called multi-static response (MSR) or scattering matrices can be handled, the total number of dipole antennas (in general, transducers and receivers) are sufficiently large enough, and measurement data is not affected by the existence of neighboring antennas. Unfortunately, in many cases, manufacturing microwave systems that can measure scattered field data with the transducer and receiver at the same location (diagonal elements) is inconvenient. Hence, designing an alternative imaging algorithm without consideration of diagonal elements of MSR or scattering matrices is an important and interesting direction for real-world application.

In this research, we apply MUSIC as a non-iterative microwave imaging algorithm when the diagonal elements of the scattering matrix are unknown. In order to show the feasibility of MUSIC in real-world microwave imaging, we carefully explore mathematical structure of imaging function by establishing a relationship with an infinite series of Bessel function of integer order, the total number and location of dipole antennas, and the applied frequency. This is based on the Born approximation in the existence of a small anomaly \cite{AK2} and the physical factorization of the scattering matrix in the presence of an extended anomaly \cite{HSZ1}. Based on the explored structure of the imaging function, we can discover various intrinsic properties and the fundamental limitations of MUSIC in real-world microwave imaging. In order to support our theoretical results, simulation results for small and various extended anomalies through the real-data generated by commercial tool and manufactured microwave machine are exhibited.

This research is organized as follows: In Section \ref{sec:2}, we briefly introduce the forward problem including scattering parameter and MUSIC algorithm without diagonal elements of scattering matrix. In Section \ref{sec:3}, we analyze mathematical structure of imaging function by establishing infinite series of Bessel functions of integer order and discover intrinsic properties of MUSIC in real-world application. In Sections \ref{sec:4} and \ref{sec:5}, experimental results respectively using synthetic and real data are exhibited to show the feasibility and limitation of the MUSIC. A conclusion including an outline of future work is given in Section \ref{sec:6}.

Finally, let us emphasize that obtained imaging results via MUSIC seems good but does not guarantee the complete shape of anomalies. Thus, one must apply iteration scheme to retrieve more accurate shape of anomaly (if there remains a discrepancy between the measurement and computed data) once an appropriate cost functional is chosen.

\section{Forward Problem and MUSIC Algorithm}\label{sec:2}
\subsection{Forward Problem and Scattering Parameter}
In this section, we briefly introduce the scattering parameter and MUSIC algorithm for imaging unknown anomalies from scattering matrix. For the sake of simplicity, we assume that there exists an anomaly $\Gamma$ with smooth boundary $\p\Gamma\in C^2$. Throughout this paper, we assume that $\Gamma$ is a small anomaly located at $\mr_\star$ and surrounded by the circular array of dipole antennas located at $\ma_n$, $n=1,2,\cdots,\N$. Throughout this paper, every targets are classified by the values of their dielectric permittivity and electrical conductivity at a given angular frequency $\omega=2\pi f$, where $f=c/\lambda$ is fixed frequency in microwave band, $c$ is the speed of light, and $\lambda$ is the wavelength. Hence, we set the value of magnetic permeability of $\Gamma$ and background are same and constant at every location $\mr$, i.e., we let $\mu(\mr)=\mu_{\mathrm{b}}=4\pi\times \SI{e-7}{\henry/\meter}$. {We denote $\eps_\star=\eps_{\rstar}\cdot\eps_0$ and $\epsb=\eps_{\rback}\cdot\eps_0$ as the permittivity of $\Gamma$ and the background, respectively and assume that assume that $\epsb\gg\sigmab/\omega$. Here, $\eps_0=\SI{8.854e-12}{\farad/\meter}$ is the vacuum permittivity.} Conductivities $\sigma_\star$ and $\sigmab$ can be defined analogously. With this, we define the following piecewise permittivity and conductivity as $\eps(\mr)$ and $\sigma(\mr)$, respectively such that
\[\eps(\mr)=\left\{\begin{array}{rcl}
                     \eps_\star & \mbox{for} & \mr\in\Gamma,\\
                     \eps_{\mathrm{b}} & \mbox{for} & \mr\in\mathbb{R}^3\backslash\overline{\Gamma},
                   \end{array}
\right.
\quad\mbox{and}\quad
\sigma(\mr)=\left\{\begin{array}{rcl}
                     \sigma_\star & \mbox{for} & \mr\in\Gamma,\\
                     \sigma_{\mathrm{b}} & \mbox{for} & \mr\in\mathbb{R}^3\backslash\overline{\Gamma},
                   \end{array}
\right.\]
With this, we introduce the background wavenumber $k$ as
\[k^2=\omega^2\mu_{\mathrm{b}}\left(\eps_{\mathrm{b}}+i\frac{\sigma_{\mathrm{b}}}{\omega}\right).\]
Let $\mE_{\inc}(\ma_n,\mr)\in\mathbb{C}^{1\times3}$ be the incident electric field in a homogeneous medium due to the point current density $\mathbf{J}$ at $\ma_n$. Then, based on the Maxwell's equation, $\mE_{\inc}(\ma_n,\mr)$ satisfies
\[\nabla\times\mE_{\inc}(\ma_n,\mr)=-i\omega\mu_{\mathrm{b}}{\mH_{\inc}}(\ma_n,\mr)\quad\mbox{and}\quad\nabla\times{\mH_{\inc}}(\ma_n,\mr)=(\sigma_{\mathrm{b}}+i\omega\eps_{\mathrm{b}})\mE_{\inc}(\ma_n,\mr),\]
where ${\mH_{\inc}}\in\mathbb{C}^{1\times3}$ denotes the magnetic field. Analogously, let $\mE_{\tot}(\mr,\ma_n)\in\mathbb{C}^{1\times3}$ be the total electric field in the existence of $\Gamma$. Then, $\mE_{\tot}(\mr,\ma_n)$ satisfies
\[\nabla\times\mE_{\tot}(\mr,\ma_n)=-i\omega\mu_{\mathrm{b}}{\mH_{\tot}}(\mr,\ma_n)\quad\mbox{and}\quad\nabla\times{\mH_{\tot}}(\mr,\ma_n)=(\sigma(\mr)+i\omega\eps(\mr))\mE_{\tot}(\mr,\ma_n)\]
with the transmission condition on $\p\Gamma$.

Let $\mathrm{S}(m,n)$ be the scattering parameter (or $S-$parameter), which is the measurement data through the microwave machine that is defined as
\[\mathrm{S}(m,n):=\frac{\mathrm{V}_m^-}{\mathrm{V}_n^+},\]
where
$\mathrm{V}_m^-$ and $\mathrm{V}_n^+$ denote the output voltage (or reflected wave) at the $m-$th antenna and the input voltage (or incident wave) at $n-$th antenna, respectively. We let $\mathrm{S}_{\tot}(m,n)$ and $\mathrm{S}_{\inc}(m,n)$ be measured scattering parameters with and without $\Gamma$, respectively. Then, the scattered field $S-$parameter $\mathrm{S}_{\scat}(m,n):=\mathrm{S}_{\tot}(m,n)-\mathrm{S}_{\inc}(m,n)$, which is the measurement data in this research\footnote{This subtraction is essential and useful to remove unknown modeling error.}, can be represented as follows (see \cite{HSM2}):
\begin{equation}\label{Sparameter}
\mathrm{S}_{\scat}(m,n)=\frac{ik^2}{4\omega\mu_{\mathrm{b}}}\int_{\Omega}\left(\frac{\eps(\mr)-\eps_{\mathrm{b}}}{\eps_{\mathrm{b}}}+i\frac{\sigma(\mr)-\sigma_{\mathrm{b}}}{\omega\eps_{\mathrm{b}}}\right)\mE_{\inc}(\ma_n,\mr)\cdot\mE_{\tot}(\mr,\ma_m)d\mr.
\end{equation}

\subsection{Introduction to MUSIC Algorithm}
Now, we introduce traditional MUSIC algorithm. For a detail, we refer to \cite{PKLS}. Generally, one assumed that the following scattering matrix can be available
\begin{equation}\label{MSR}
\mathbb{K}=\left[
               \begin{array}{ccccc}
                  \mathrm{S}_{\scat}(1,1) & \mathrm{S}_{\scat}(1,2) & \cdots & \mathrm{S}_{\scat}(1,\N-1) & \mathrm{S}_{\scat}(1,\N) \\
                  \mathrm{S}_{\scat}(2,1) & \mathrm{S}_{\scat}(2,2) & \cdots & \mathrm{S}_{\scat}(2,\N-1) & \mathrm{S}_{\scat}(2,\N) \\
                  \vdots & \vdots & \ddots & \vdots & \vdots \\
                 \mathrm{S}_{\scat}(\N,1) & \mathrm{S}_{\scat}(\N,2) & \cdots & \mathrm{S}_{\scat}(\N,\N-1) & \mathrm{S}_{\scat}(\N,\N)
               \end{array}
             \right].
\end{equation}
{This means that the entries of $\mathbb{K}$ are obtained while the other ports are matched.} The main purpose of MUSIC is to find locations of $\mr_\star\in\Gamma$ from $\mathbb{K}$ without any \textit{a priori} information of $\Gamma$. Based on the representation \eqref{Sparameter}, we need the value of total field $\mE_{\tot}(\mr,\ma_m)$ at every $\mr$. Unfortunately, direct calculation of total field is impossible unless we have complete information of $\Gamma$. Thus, to design the imaging algorithm of MUSIC, we take the next two steps to consider an alternative representation of \eqref{Sparameter}.

\begin{rem}[Linearization]
Now, let us assume that
\begin{equation}\label{SmallAnomaly}
\left(\sqrt{\frac{\eps_\star}{\eps_{\mathrm{b}}}}-1\right)r<\frac{\lambda}{4}.
\end{equation}
Then, based on \cite{SKL}, $\Gamma$ can be regarded as a small anomaly so that by applying the Born approximation \cite{HSM2} and reciprocity property of incident field, (\ref{Sparameter}) can be approximated as follows:
\begin{align}
\begin{aligned}\label{Formula-Sparameter}
\mathrm{S}_{\scat}(m,n)&=\frac{ik^2}{4\omega\mu_{\mathrm{b}}}\int_{\Omega}\left(\frac{\eps(\mr)-\eps_{\mathrm{b}}}{\eps_{\mathrm{b}}}+i\frac{\sigma(\mr)-\sigma_{\mathrm{b}}}{\omega\eps_{\mathrm{b}}}\right)\mE_{\inc}(\ma_n,\mr)\cdot\mE_{\tot}(\mr,\ma_m)d\mr\\
&\approx\frac{ik^2}{4\omega\mu_{\mathrm{b}}}\int_{\Gamma}\left(\frac{\eps(\mr)-\eps_{\mathrm{b}}}{\eps_{\mathrm{b}}}+i\frac{\sigma(\mr)-\sigma_{\mathrm{b}}}{\omega\eps_{\mathrm{b}}}\right)\mE_{\inc}(\ma_n,\mr)\cdot\mE_{\inc}(\ma_m,\mr)d\mr.
\end{aligned}
\end{align}
\end{rem}

\begin{rem}[2D approximation]
Notice that the height of microwave machine can be said to be long enough (see \cite{P-SUB11,SLP} and Figure \ref{Machine} for experimental setup). Hence, based on the mathematical treatment of the scattering of time-harmonic electromagnetic waves from thin infinitely long cylindrical obstacles (see \cite{K} for instance), only $z-$component of $\mE_{\inc}(\ma_n,\mr)$ and $\mE_{\tot}(\mr,\ma_m)$ can be handled for $n,m=1,2,\cdots,\N$, and this problem becomes a two-dimensional problem, i.e., we can image the cross-section of $\Gamma$. Now, let us denote $\E_{\inc}^{(z)}(\ma_n,\mr)$ and $\E_{\tot}^{(z)}(\mr,\ma_m)$ be the $z-$component of incident and total fields, respectively. Then, \eqref{Formula-Sparameter} can be written as
\begin{align}
\begin{aligned}\label{Sparameter-Representation}
\mathrm{S}_{\scat}(m,n)&\approx\frac{ik^2}{4\omega\mu_{\mathrm{b}}}\int_{\Gamma}\left(\frac{\eps(\mr)-\eps_{\mathrm{b}}}{\eps_{\mathrm{b}}}+i\frac{\sigma(\mr)-\sigma_{\mathrm{b}}}{\omega\eps_{\mathrm{b}}}\right)\mE_{\inc}(\ma_n,\mr)\cdot\mE_{\inc}(\ma_m,\mr)d\mr\\
&\approx\frac{ik^2\area(\Gamma)}{4\omega\mu_{\mathrm{b}}}\left(\frac{\eps_\star-\eps_{\mathrm{b}}}{\eps_{\mathrm{b}}}+i\frac{\sigma_\star-\sigma_{\mathrm{b}}}{\omega\eps_{\mathrm{b}}}\right)\E_{\inc}^{(z)}(\ma_n,\mr_\star)\E_{\inc}^{(z)}(\ma_m,\mr_\star).
\end{aligned}
\end{align}
Here, $\E_{\inc}^{(z)}(\ma_n,\mr)$ is given by
\[\E_{\inc}^{(z)}(\ma_n,\mr)=-\frac{i}{4}H_0^{(1)}(k|\ma_n-\mr|),\]
where $H_0^{(1)}$ denotes the Hankel function of order zero of the first kind (see \cite{P-SUB11,SLP} for instance).
\end{rem}

Based on the representation \eqref{Sparameter-Representation}, we can examine that the range of scattering matrix $\mathbb{K}$ can be determined from the span of
\begin{equation}\label{SpanVectors}
\mW(\mr_\star)=\bigg[\E_{\inc}^{(z)}(\ma_1,\mr_\star),\E_{\inc}^{(z)}(\ma_2,\mr_\star),\cdots,\E_{\inc}^{(z)}(\ma_{\N},\mr_\star)\bigg]^{\mathtt{T}}
\end{equation}
corresponding to $\Gamma$. To introduce the MUSIC algorithm, we apply singular value decomposition (SVD) to $\mathbb{K}$, as follows:
\[\mathbb{K}=\sum_{n=1}^{\N}\tau_n\mU_n\mV_n^*\approx\sum_{n=1}^{\M}\tau_n\mU_n\mV_n^*,\]
where $\tau_n$ are the singular values that satisfy
\[\tau_1\geq\tau_2\geq\cdots\geq\tau_{\M}>0\quad\text{and}\quad\tau_n\approx0\quad\text{for}\quad n=\M+1,\M+2,\cdots,\N,\]
and $\mU_n$ and $\mV_n$ are respectively the left- and right-singular vectors of $\mathbb{K}$. Then, the first $\M$ columns of left-singular vectors $\set{\mU_1,\mU_2,\cdots,\mU_{\M}}$ provide an orthonormal basis for $\mathbb{K}$. With this, we introduce a projection operator onto the noise subspace as
\begin{equation}\label{Projection}
\mathbb{P}_{\noise}=\mathbb{I}_{\N}-\sum_{n=1}^{\M}\mU_n\mU_n^*,
\end{equation}
where $\mathbb{I}_{\N}$ denotes $\N\times\N$ identity matrix. Now, let us consider the normalized vector $\mf(\mr)=\mW(\mr)/|\mW(\mr)|$, where $\mW$ is given by \eqref{SpanVectors}. Then, based on \cite{AK2,K2}, the following statement holds:
\[{\mf(\mr)\in\mbox{Range}(\mathbb{K})\quad\mbox{if and only if}\quad \mr=\mr_\star\in\Gamma.}\]
This means that, $|\mathbb{P}_{\noise}(\mf(\mr))|=0$ when $\mr\in\Gamma$. Hence, the location of $\Gamma$ can be identified by plotting the traditional MUSIC-type imaging function
\begin{equation}\label{TraditionalImagingFunctionMUSIC}
  \mathfrak{F}_{\tm}(\mr)=\frac{1}{|\mathbb{P}_{\noise}(\mf(\mr))|},\quad\mr\in\Omega.
\end{equation}
The plot of $\mathfrak{F}_{\tm}(\mr)$ is expected to exhibit peaks of large magnitude (in theory $+\infty$) at $\mr=\mr_\star\in\Gamma$.

\begin{rem}[Total number of nonzero singular values]
Theoretically, total number of nonzero singular values $\M$ of $\mathbb{K}$ must be equal to the total number of anomalies $($if there exists a small anomaly, $\M=1)$. However, in real-world application, $\M>1$ because the measurement data $\mathrm{S}_{\scat}(m,n)$ is influenced by not only anomaly but also antennas. Sometimes, this phenomenon is called as the ``multiple scattering effect''. We refer to \cite{P-SUB11} for a detailed discussion. {It is worth mentioning that for 3D small targets for each targets there are as much as three singular values.}
\end{rem}

\begin{rem}[Extended anomaly]
If the condition \eqref{SmallAnomaly} does not satisfied, $\Gamma$ is regarded as an extended anomaly. Since $\Gamma$ is not small, we cannot apply the Born approximation so that we cannot design MUSIC algorithm with the similar manner. In this case, based on the physical factorization of the scattering matrix \cite{HSZ1}, $\mathbb{K}$ can be decomposed as
\begin{align}
\begin{aligned}\label{Factorization}
  \mathbb{K}&=\frac{ik^2}{4\omega\mu_{\mathrm{b}}}\int_{\Omega}\left(\frac{\eps(\mr)-\eps_{\mathrm{b}}}{\eps_{\mathrm{b}}}+i\frac{\sigma(\mr)-\sigma_{\mathrm{b}}}{\omega\sigma_{\mathrm{b}}}\right)\mathbb{E}_{\inc}(\mr)\mathbb{E}_{\tot}(\mr)d\mr\\
  &=\frac{ik^2}{4\omega\mu_{\mathrm{b}}}\int_{\Gamma}\left(\frac{\eps_\star-\eps_{\mathrm{b}}}{\eps_{\mathrm{b}}}+i\frac{\sigma_\star-\sigma_{\mathrm{b}}}{\omega\eps_{\mathrm{b}}}\right)\mathbb{E}_{\inc}(\mr)\mathbb{E}_{\tot}(\mr)d\mr,
  \end{aligned}
  \end{align}
  where
  \[\mathbb{E}_{\inc}(\mr)=\bigg[\E_{\inc}^{(z)}(\ma_1,\mr),\E_{\inc}^{(z)}(\ma_2,\mr),\cdots,\E_{\inc}^{(z)}(\ma_\N,\mr)\bigg]^{\mathtt{T}}\quad\text{and}\quad\mathbb{E}_{\tot}(\mr)=\bigg[\E_{\tot}^{(z)}(\ma_1,\mr),\E_{\tot}^{(z)}(\ma_2,\mr),\cdots,\E_{\tot}^{(z)}(\ma_\N,\mr)\bigg],\]
respectively. Based on this factorization, the range of $\mathbb{K}$ is determined based on the span of $\mathbb{E}_{\inc}(\mr)$ corresponding to $\mr\in\Gamma$. This means that the imaging function $\mathfrak{F}_{\tm}(\mr)$ of (\ref{TraditionalImagingFunctionMUSIC}) can be defined by selecting the first $\M-$singular vectors of $\mathbb{K}$, and the outline shape of $\Gamma$ will be imaged through the $\mathfrak{F}_{\tm}(\mr)$. We refer to \cite{HSZ1} for instance and Figure \ref{CST3} for a related result.
\end{rem}

\subsection{MUSIC Algorithm Without Diagonal Elements of the Scattering Matrix}
Unlike the conventional researches, we cannot measure $\mathrm{S}_{\scat}(n,n)$ for $n=1,2,\cdots,\N$, because each of the $\N$ antennas is used for signal transmission and the remaining $\N-1$ antennas are used for signal reception (see \cite{SSKLJ} for instance) so that in real-world problem, obtained scattering matrix is given by
\begin{equation}\label{MSRunknown}
\mathbb{S}=\left[
               \begin{array}{ccccc}
                  \mbox{unknown} & \mathrm{S}_{\scat}(1,2) & \cdots & \mathrm{S}_{\scat}(1,\N-1) & \mathrm{S}_{\scat}(1,\N) \\
                  \mathrm{S}_{\scat}(2,1) & \mbox{unknown} & \cdots & \mathrm{S}_{\scat}(2,\N-1) & \mathrm{S}_{\scat}(2,\N) \\
                  \vdots & \vdots & \ddots & \vdots & \vdots \\
                 \mathrm{S}_{\scat}(\N,1) & \mathrm{S}_{\scat}(\N,2) & \cdots & \mathrm{S}_{\scat}(\N,\N-1) & \mbox{unknown}
               \end{array}
             \right].
\end{equation}
In this case, it is impossible to determine the range of scattering matrix \eqref{MSRunknown} so that we cannot apply traditional MUSIC algorithm for imaging $\Gamma$. Motivated by \cite{P-SUB11}, we simply set $\mathrm{S}_{\scat}(n,n)=0$ for all $n$ in the scattering matrix \eqref{MSRunknown} and apply MUSIC algorithm, i.e., we consider the following matrix:
\begin{equation}\label{MSRZERO}
\mathbb{D}=\left[
               \begin{array}{ccccc}
                  0 & \mathrm{S}_{\scat}(1,2) & \cdots & \mathrm{S}_{\scat}(1,\N-1) & \mathrm{S}_{\scat}(1,\N) \\
                  \mathrm{S}_{\scat}(2,1) & 0 & \cdots & \mathrm{S}_{\scat}(2,\N-1) & \mathrm{S}_{\scat}(2,\N) \\
                  \vdots & \vdots & \ddots & \vdots & \vdots \\
                 \mathrm{S}_{\scat}(\N,1) & \mathrm{S}_{\scat}(\N,2) & \cdots & \mathrm{S}_{\scat}(\N,\N-1) & 0
               \end{array}
             \right].
\end{equation}
In this case, determination of the range of \eqref{MSRZERO} is still impossible. Nevertheless, we can identify an outline shape of $\Gamma$ through the MUSIC algorithm: we apply singular value decomposition (SVD) to $\mathbb{D}$, as follows:
\[\mathbb{D}=\sum_{n=1}^{\N}\tau_n\mU_n\mV_n^*\approx\tau_1\mU_1\mV_1^*\]
because there is no multiple scattering effect (see \cite{P-SUB11} for instance). With this, we introduce a projection operator onto the noise subspace as
\[\mathbb{P}_{\noise}=\mathbb{I}_{\N}-\mU_1\mU_1^*,\]
and corresponding MUSIC-type imaging function without diagonal elements of $\mathbb{D}$
\begin{equation}\label{ImagingFunctionMUSIC}
  \mathfrak{F}_{\dm}(\mr)=\frac{1}{|\mathbb{P}_{\noise}(\mf(\mr))|}.
\end{equation}
Then, it is possible to identify the outline shape of $\Gamma$ through the map of $\mathfrak{F}_{\dm}(\mr)$. Theoretical reason of applicability including unique determination is discussed in the next section.

\section{Analysis of the Imaging Function}\label{sec:3}
Here, we carefully derive mathematical structure of $\mathfrak{F}_{\dm}(\mr)$ by establishing a relationship with an infinite series of Bessel function of integer order and examine intrinsic properties of MUSIC. For this purpose, we carefully explore mathematical structure of $\mathfrak{F}_{\dm}(\mr)$. The derivation is given in Section \ref{sec:A}.
\begin{thm}[Structure of imaging function for single anomaly]\label{StructureImaging}
Let $\vt_n=\ma_n/|\ma_n|=\ma_n/R=[\cos\theta_n,\sin\theta_n]^{\mathtt{T}}$ and $\mr-\mr_\star=|\mr-\mr_\star|[\cos\phi_\star,\sin\phi_\star]^{\mathtt{T}}$. {If $\ma_n$ satisfies $|\ma_n-\mr|\gg1/4|k|$ for $n=1,2,\cdots,\N$,} $\mathfrak{F}_{\dm}(\mr)$ can be represented as follows:
\begin{equation}\label{StructureImagingFunction}
\mathfrak{F}_{\dm}(\mr)\approx\frac{\N-1}{\N}\abs{1-J_0(k|\mr_\star-\mr|)^2-\frac{1}{\N}\mathrm{Re}(\mathcal{E}(\mr,\mr_\star))}^{-1/2},
\end{equation}
where $J_\nu$ denotes the Bessel function of integer order $\nu$ of the first kind and $\mathrm{Re}(\mathcal{E}(\mr,\mr_\star))$ is the real-part of residual term $\mathcal{E}(\mr,\mr_\star)=O(1)$, which is given by
\begin{align*}
\mathcal{E}(\mr,\mr_\star)=&2J_0(k|\mr_\star-\mr|)\sum_{n=1}^{\N}\sum_{\nu\in\mathbb{Z}^*}i^\nu J_{\nu}(k|\mr_\star-\mr|)\e^{i\nu(\theta_n-\phi_\star)}\\
&+\frac{1}{\N}\left(\sum_{n=1}^{\N}\sum_{\nu\in\mathbb{Z}^*}i^\nu J_{\nu}(k|\mr_\star-\mr|)\e^{i\nu(\theta_n-\phi_\star)}\right)\left(\sum_{n=1}^{\N}\sum_{\nu\in\mathbb{Z}^*}(-i)^\nu J_{\nu}(k|\mr_\star-\mr|)\e^{-i\nu(\theta_n-\phi_\star)}\right).
\end{align*}
Here, $\mathbb{Z}$ be the set of integer number and $\mathbb{Z}^*=\mathbb{Z}\cup\set{+\infty,-\infty}\backslash\set{0}$. Furthermore, for all $n=1,2,\cdots,\N$,
\[\lim_{\mr\to\ma_n}\mathfrak{F}_{\dm}(\mr)=0.\]
\end{thm}

Based on the identified structure \eqref{StructureImagingFunction}, we can examine several properties of the imaging function $\mathfrak{F}_{\dm}(\mr)$.

\begin{rem}[Applicability]\label{Remark-Applicability}
Since $J_0(k|\mr_\star-\mr|)=1$ and $\mathcal{E}(\mr,\mr_\star)=0$ when $\mr=\mr_\star\in\Gamma$, $\mathfrak{F}_{\dm}(\mr)$ will have it maximum at $\mr_\star\in\Gamma$. Thus, it will be possible to identify location of $\Gamma$ via the map of $\mathfrak{F}_{\dm}(\mr)$.
\end{rem}

Based on the Remark \ref{Remark-Applicability}, we can obtain the following result of unique determination.
\begin{cor}[Unique determination of anomaly]\label{UniqueDetermination}
   For given angular frequency $\omega$, the location of anomaly can be identified uniquely through the map of $\mathfrak{F}_{\dm}(\mr)$.
\end{cor}

\begin{rem}[Appearance of unexpected artifacts and resolution]
If $\mr\ne\mr_\star$, the terms $J_0(k|\mr_\star-\mr|)$ and $\mathcal{E}(\mr,\mr_\star)$ are the main causative agents of generating unexpected artifacts. The resolution of imaging result is significantly depending on the value of applied frequency $\omega$. If $\omega$ is sufficiently large, i.e., if one applies high frequency, one can obtain a good imaging result with high resolution but several artifacts will be contained in the imaging result also. In contrast, if one applies low frequency, appearance of artifacts can be reduced but the resolution will be poor.
\end{rem}

\begin{rem}[Ideal condition for a better imaging performance]
 Notice that if $N\longrightarrow+\infty$ then, $\mathcal{E}(\mr,\mr_\star)/\N\longrightarrow0$ and correspondingly, $\mathfrak{F}_{\dm}(\mr)\longrightarrow+\infty$ because $|\mathbb{P}_{\noise}(\mW(\mr))|\longrightarrow0$. This means that increasing total number of antennas guarantees the better imaging performance. Unfortunately, increasing total number of antennas is difficult problem due to the mechanical circumstance. If $\omega\longrightarrow+\infty$, one will obtain good result but application of infinite-valued frequency is an ideal condition.
\end{rem}

\begin{rem}[Experimental configuration for a better imaging performance]
Based on the structure of imaging function \eqref{StructureImagingFunction}, it is easy to see that the term $\mathcal{E}(\mr,\mr_\star)$ obstructs imaging performance. Hence, eliminating this term will guarantee good result. Notice that we cannot increase $\N$ and control $J_\nu(k|\mr_\star-\mr|)=0$ for all $\mr\in\Omega$. This means that we must focus on the elimination of $\e^{i\nu(\theta_n-\phi_\star)}$ and $\e^{-i\nu(\theta_n-\phi_\star)}$. Since we have no information of location of target, the value of $\phi_\star$ is still unknown. {Fortunately, based on the following elementary conditions of trigonometric functions
\[\cos(\theta+\pi)=-\cos\theta\quad\mbox{and}\quad\sin(\theta+\pi)=-\sin\theta,\]
we can examine that if there exists two antennas located at $\ma_1=|R|[\cos\theta_1,\sin\theta_1]^{\mathtt{T}}$ and $\ma_2=|R|[\cos\theta_2,\sin\theta_2]^{\mathtt{T}}=|R|[\cos(\theta_1+\pi),\sin(\theta_1+\pi)]^{\mathtt{T}}=-\ma_1$ then
\[\sum_{n=1}^{2}\sum_{\nu\in\mathbb{Z}^*}i^\nu J_{\nu}(k|\mr_\star-\mr|)\e^{i\nu(\theta_n-\phi_\star)}=0.\]
This means that to reduce the effect of $\mathrm{Re}(\mathcal{E}(\mr,\mr_\star))$ for obtaining a good result, the total number $\N$ of the antennas must be even and they must be distributed uniformly on the circle of radius $R$ and symmetry from each others.}
\end{rem}

The mathematical setting and analysis could be extended straightforwardly to multiple anomalies. The derivation is not provided herein; only theoretical result is shown as follows. Of course, unique determination is guaranteed for multiple anomalies.

\begin{cor}[Structure of imaging function for multiple anomalies]
Assume that there exists multiple anomalies $\Gamma_s$ with location $\mr_s$, $s=1,2,\cdots,S$. Let $\vt_n=\ma_n/|\ma_n|=\ma_n/R=[\cos\theta_n,\sin\theta_n]^{\mathtt{T}}$ and $\mr-\mr_s=|\mr-\mr_s|[\cos\phi_s,\sin\phi_s]^{\mathtt{T}}$. If $\ma_n$ satisfies $|\ma_n-\mr|\gg1/4|k|$ for $n=1,2,\cdots,\N$, $\mathfrak{F}_{\dm}(\mr)$ can be represented as follows: if $\mr$ is far from $\ma_n$,
\begin{equation}\label{StructureImagingFunctionMultiple}
\mathfrak{F}_{\dm}(\mr)\approx\frac{\N-1}{\N}\abs{\sum_{s=1}^{S}\left(1-J_0(k|\mr_s-\mr|)^2-\frac{1}{\N}\mathrm{Re}(\mathcal{E}(\mr,\mr_s))\right)}^{-1/2},
\end{equation}
where $\mathcal{E}(\mr,\mr_s)$ is introduced in Theorem \ref{StructureImaging}. Furthermore, for all $n=1,2,\cdots,\N$,
\[\lim_{\mr\to\ma_n}\mathfrak{F}_{\dm}(\mr)=0.\]
\end{cor}

\section{Simulation Results With Synthetic Data: CST STUDIO SUITE}\label{sec:4}
First, we exhibit several results of numerical simulation with synthetic data. For this, we use $\N=16$ dipole antennas located at
\[\ma_n=\SI{0.09}{\meter}[\cos\theta_n,\sin\theta_n]^{\mathtt{T}},\quad\theta_n=\frac{3\pi}{2}-\frac{2\pi(n-1)}{\N},\]
and we apply $f=\SI{1.2}{\giga\hertz}$ frequency. To perform the imaging, we set the ROI $\Omega$ to $\Omega=\SI{0.1}{\meter}[-1,1]^{\mathtt{T}}\times\SI{0.1}{\meter}[-1,1]^{\mathtt{T}}$.

The values of permittivity, conductivity, location, and size of each circle-like anomaly are given in Table \ref{table}. With this, by checking \eqref{SmallAnomaly}, $\Gamma_\mathrm{S}$ and $\Gamma_\mathrm{E}$ in Table \ref{table} are regarded as a small and extended anomaly, respectively, refer to \cite{SKL,P-SUB11}.

\begin{table}[h]
\begin{center}
\begin{tabular}{c||c|c|c|c}
\hline\hline  \centering Target&Permittivity&Conductivity ($\SI{}{\siemens/\meter}$)&Location&Radius\\
\hline\hline\centering Background&$20$ ($\eps_{\rback}$)&$0.2$ ($\sigmab$)&$-$&$-$\\
\hline\centering Small anomaly $\Gamma_\mathrm{S}$&$55$ ($\eps_{\rstar}$)&$1.2$ ($\sigma_{\star}$)&$\SI{0.01}{\meter}[1,3]^{\mathtt{T}}$&$\SI{0.01}{\meter}$\\
\hline\centering Extended anomaly $\Gamma_\mathrm{E}$&$15$ ($\eps_{\rstar}$)&$0.5$ ($\sigma_{\star}$)&$\SI{0.01}{\meter}[1,2]^{\mathtt{T}}$&$\SI{0.05}{\meter}$\\
\hline\hline
\end{tabular}
\caption{\label{table}Values of permittivities, conductivities, locations, and sizes of anomalies as well as the background.}
\end{center}
\end{table}

\begin{figure}[h]
  \centering
  \includegraphics[width=0.495\textwidth]{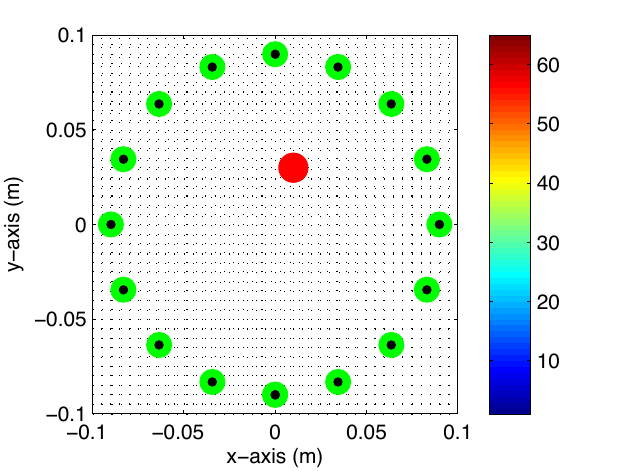}
  \includegraphics[width=0.495\textwidth]{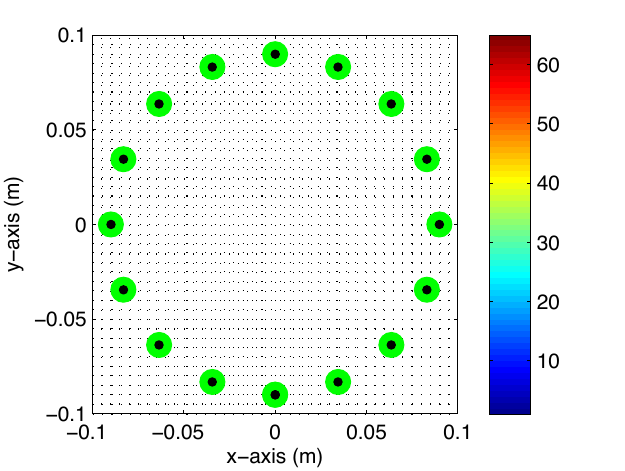}
  \caption{\label{Setting-Small}(Examples \ref{Example4-1} and \ref{Example4-2}) Test configurations with (left) and without (right) anomaly $\Gamma_\mathrm{S}$ for calculating $\mathrm{S}_{\tot}(m,n)$ and $\mathrm{S}_{\inc}(m,n)$, respectively.}
\end{figure}

\begin{figure}[h]
  \centering
  \includegraphics[width=0.495\textwidth]{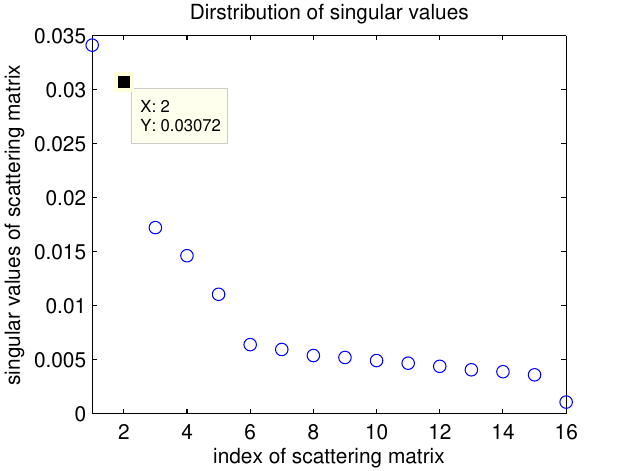}
  \includegraphics[width=0.495\textwidth]{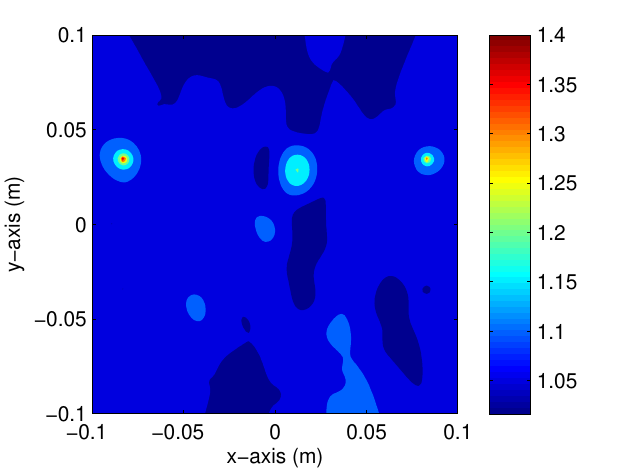}\\
  \includegraphics[width=0.495\textwidth]{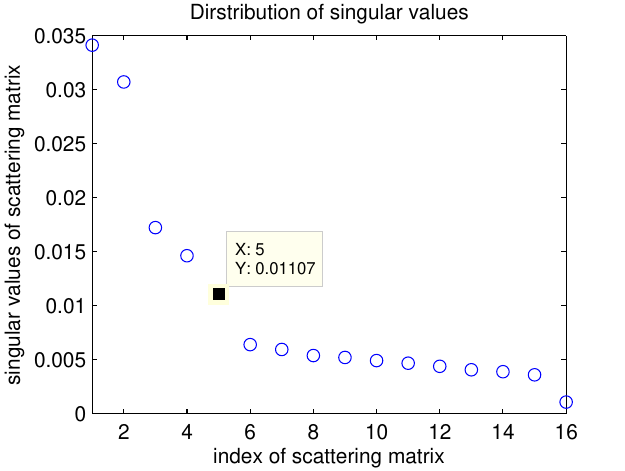}
  \includegraphics[width=0.495\textwidth]{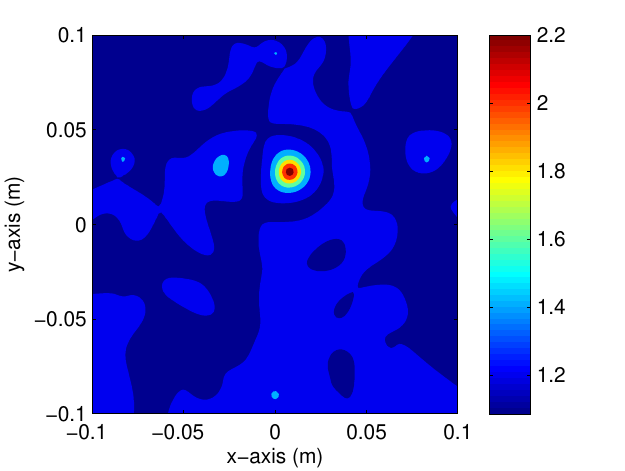}
  \caption{\label{CST1}(Example \ref{Example4-1}) Distribution of singular values of scattering matrix $\mathbb{K}$ (left column) and map of $\mathfrak{F}_{\tm}(\mr)$ (right column).}
\end{figure}

\begin{ex}[Imaging of a small anomaly: selection of nonzero singular values]\label{Example4-1}
First, let us consider the effect on the discrimination of nonzero singular values of scattering matrix $\mathbb{K}$ of \eqref{MSR}. We refer to Figure \ref{Setting-Small} for illustrating simulation setup. Figure \ref{CST1} shows distribution of singular values of $\mathbb{K}$ and maps of $\mathfrak{F}_{\tm}(\mr)$ with $\M=2$ and $\M=5$ when the anomaly is $\Gamma_\mathrm{S}$. It is well-known that total number of nonzero singular values must be equal to $1$ when there exists small anomaly. However, opposite to the traditional results, total number of nonzero singular values is not equal to the number of small anomalies anymore. Based on \cite{P-SUB11}, diagonal elements were affected by the anomaly $\Gamma_\mathrm{S}$ and antennas $\ma_n$ so that peaks of large magnitudes are shown in the map of $\mathfrak{F}_{\tm}(\mr)$ with $\M=2$ and they disturb the recognition of location of $\Gamma_\mathrm{S}$. Selection $\M=5$ yields good result but since we have \textit{a priori} information of target, nonzero singular values must be chosen very carefully when diagonal elements of scattering matrix are existing.
\end{ex}

\begin{figure}[h]
  \centering
  \includegraphics[width=0.495\textwidth]{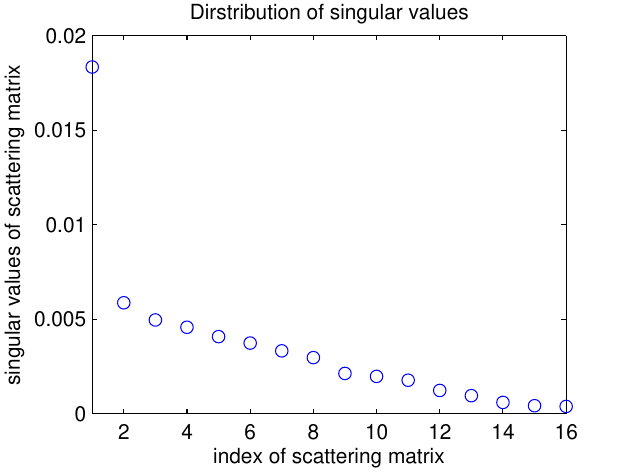}
  \includegraphics[width=0.495\textwidth]{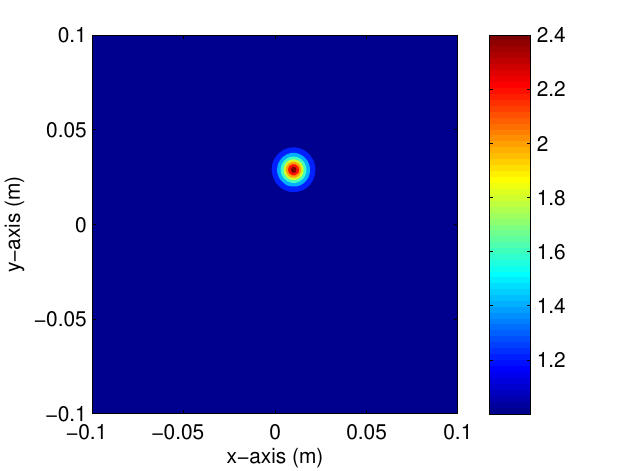}
  \caption{\label{CST2}(Example \ref{Example4-2}) Distribution of singular values of scattering matrix $\mathbb{D}$ (left) and map of $\mathfrak{F}_{\dm}(\mr)$ (right).}
\end{figure}

\begin{ex}[Imaging of a small anomaly: effects on the diagonal elements]\label{Example4-2}
Now, let us consider the imaging of $\Gamma_\mathrm{S}$ via  $\mathbb{D}$ of \eqref{MSRZERO}. Based on the result in Figure \ref{CST2}, same as the theory, it is possible to observe that total number of nonzero singular values equal to the number of anomaly because every elements of $\mathbb{D}$ are influenced by the anomaly $\Gamma_\mathrm{S}$ only. Furthermore, by comparing results in Figure \ref{CST1}, map of $\mathfrak{F}_{\dm}(\mr)$ via $\mathbb{D}$ yields good results because unexpected artifacts were disappeared. Hence, in contrast to the traditional results, eliminating diagonal elements of scattering matrix will guarantee good results.
\end{ex}

\begin{figure}[h]
  \centering
  \includegraphics[width=0.495\textwidth]{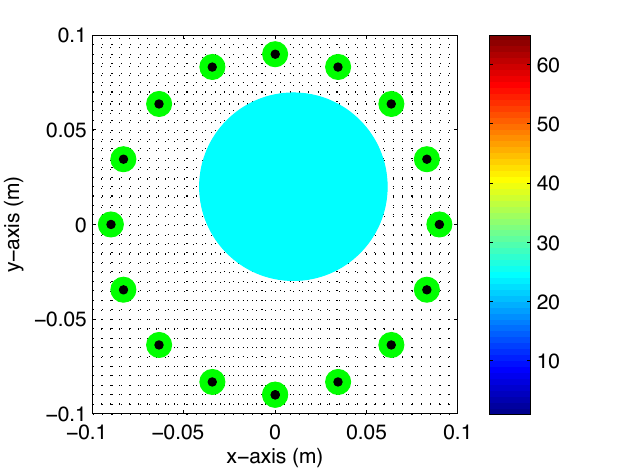}
  \includegraphics[width=0.495\textwidth]{Configuration-homo}
  \caption{\label{Setting-Extended}(Example \ref{Example4-3}) Same as Figure \ref{Setting-Small} except the anomaly is $\Gamma_\mathrm{E}$.}
\end{figure}

\begin{figure}[h]
  \centering
  \includegraphics[width=0.495\textwidth]{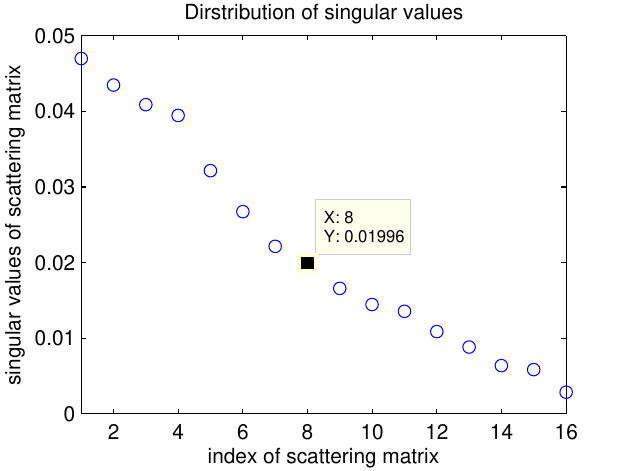}
  \includegraphics[width=0.495\textwidth]{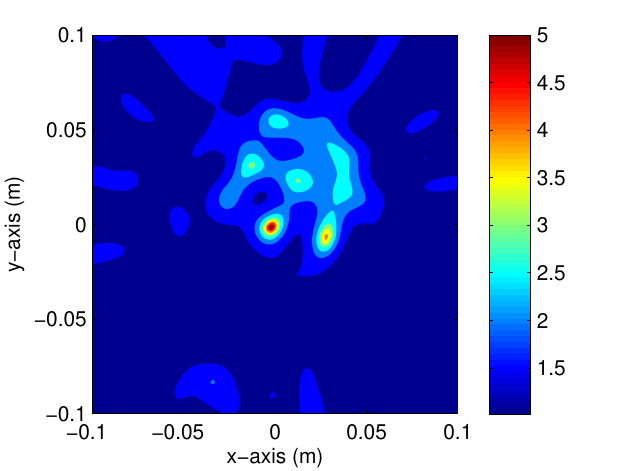}\\
  \includegraphics[width=0.495\textwidth]{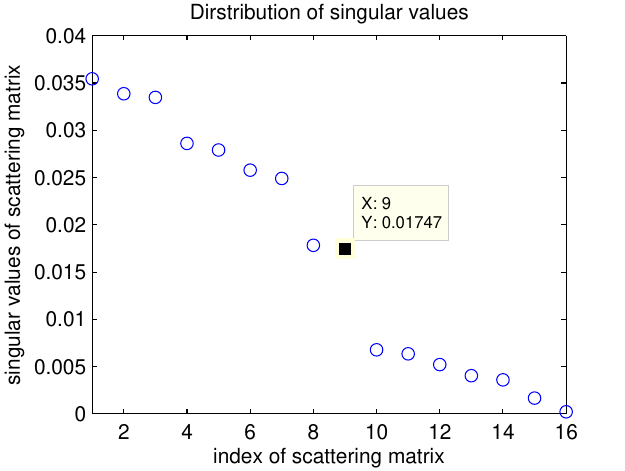}
  \includegraphics[width=0.495\textwidth]{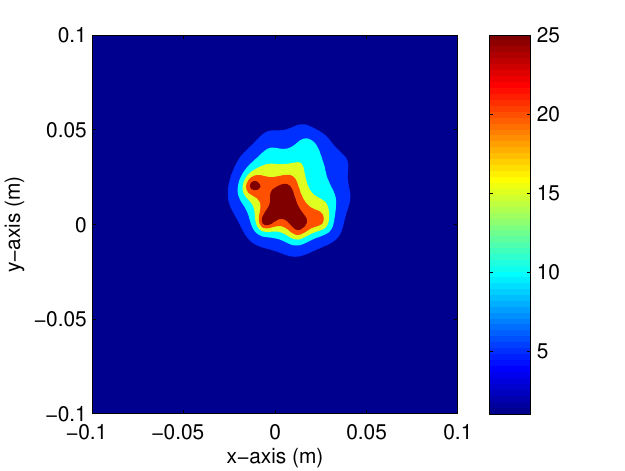}
  \caption{\label{CST3}(Example \ref{Example4-3}) Distribution of singular values of scattering matrix (left column) and map of $\mathfrak{F}_{\tm}(\mr)$ (top, right column) and $\mathfrak{F}_{\dm}(\mr)$ (bottom, right).}
\end{figure}

\begin{ex}[Imaging of a large anomaly]\label{Example4-3}
Now, we consider the imaging of extended anomaly, refer to Figure \ref{Setting-Extended} for simulation configuration. {Notice that for imaging of an extended target, the selection of nonzero singular values is very important and challenging issue, refer to \cite{HSZ1}. Fortunately, we can select $\M=9$ singular values from $\mathbb{D}$ for defining the projection operator but it is hard to select nonzero singular values from $\mathbb{K}$. In this Example, similar to the recent work \cite{PKLS}, we applied the $0.2-$thresholding scheme (choosing first $m$ singular values $\tau_n$ such that $\tau_m/\tau_1\geq0.4$) for selecting singular values of $\mathbb{K}$.} Based on the imaging results $\mathfrak{F}_{\tm}(\mr)$ and $\mathfrak{F}_{\dm}(\mr)$ for $\Gamma_\mathrm{E}$ in Figure \ref{CST3}, we can observe that although it is impossible to obtain the complete shape of an extended anomaly, an outline shape of $\Gamma_\mathrm{E}$ can be recognized via the map of $\mathfrak{F}_{\dm}(\mr)$ with $\mathbb{D}$ while it is impossible to obtain any information of $\Gamma_\mathrm{E}$ via the map of $\mathfrak{F}_{\tm}(\mr)$ with $\mathbb{K}$. Hence, we can conclude that eliminating diagonal elements of scattering matrix is a method of improvement over the traditional MUSIC algorithm.
\end{ex}

\section{Real-Data Experiments}\label{sec:5}
\begin{figure}[h]
  \centering
  \includegraphics[width=0.155\textwidth]{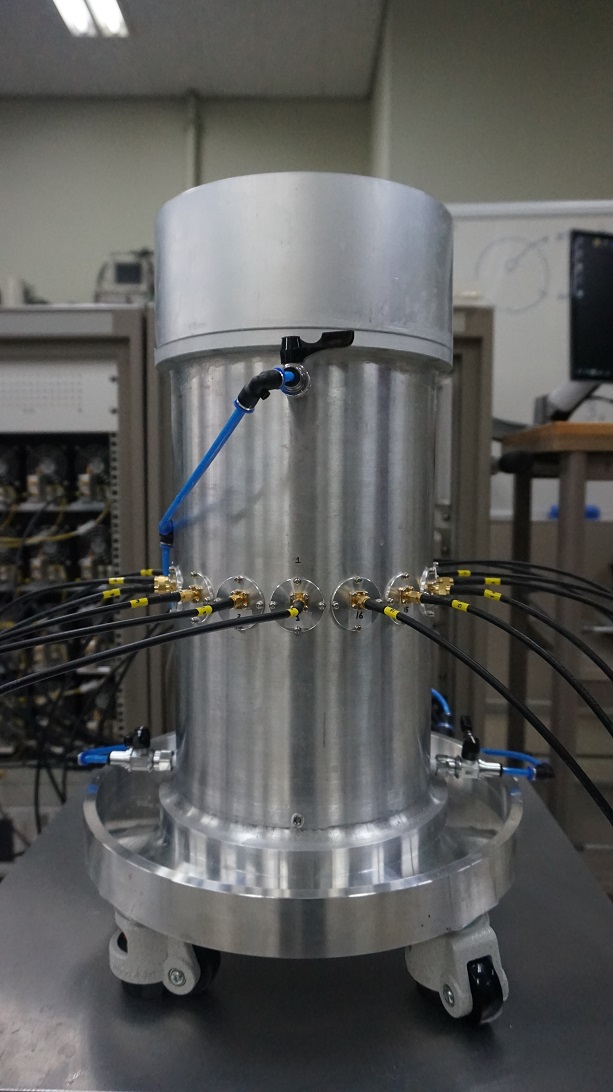}
  \includegraphics[width=0.495\textwidth]{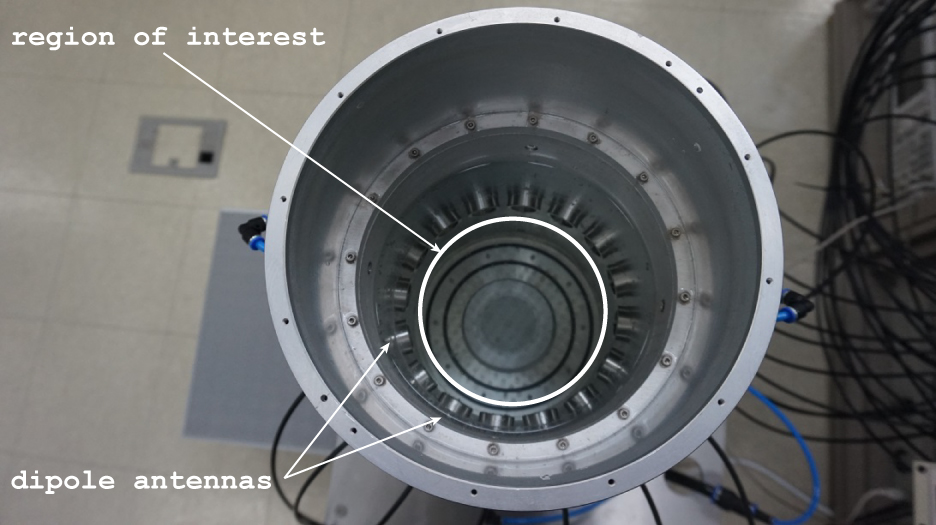}
  \caption{\label{Machine}Developed microwave machine and region of interests (ROI).}
\end{figure}

Now, we exhibit imaging results with real date to demonstrate the feasibility and support theoretical result. In order to generate $\mathrm{S}_{\tot}(m,n)$ and $\mathrm{S}_{\inc}(m,n)$, the microwave machine manufactured by the research team of the Radio Environment \& Monitoring research group of the Electronics and Telecommunications Research Institute (ETRI) was used. For a detailed description of the microwave machine, we refer to \cite{KLKJS}. For the simulation, we filled this machine by matching liquid (water) with permittivity $\eps_{\rback}=78$ and conductivity $\sigmab=\SI{0.2}{\siemens/\meter}$, apply $f=\SI{925}{\mega\hertz}$ frequency, and select the ROI as a circle of radius $\SI{0.085}{\meter}$ {to satisfy the condition in Theorem \ref{StructureImaging}}, refer to Figure \ref{Machine}. It is worth mentioning again that the scattering matrix $\mathbb{K}$ of \eqref{MSR} cannot be available. 

\begin{figure}[h]
  \centering
  \includegraphics[width=0.495\textwidth]{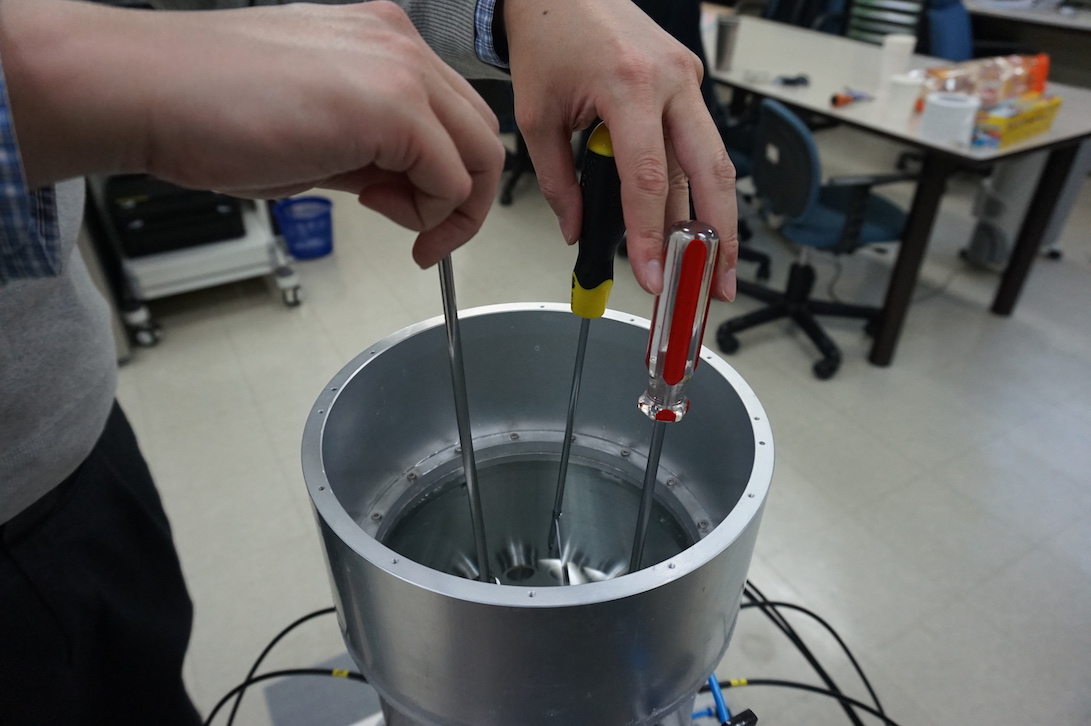}
  \includegraphics[width=0.495\textwidth]{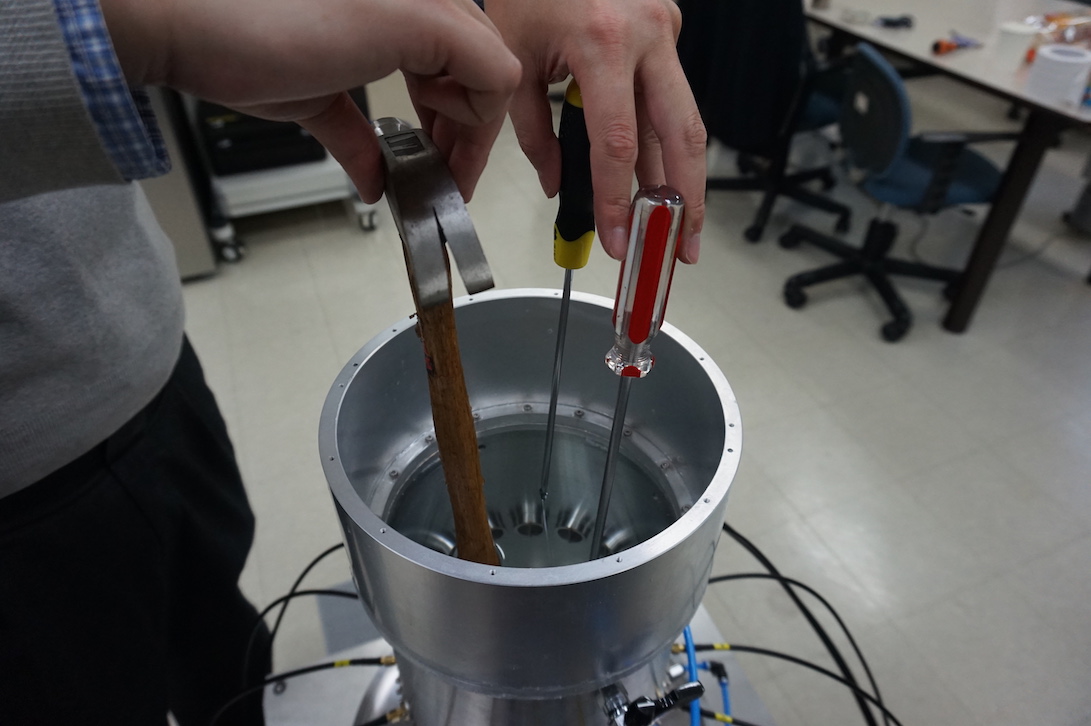}
  \caption{\label{Configuration-Real1}Real-data experiments for Examples \ref{Example5-1} (left) and \ref{Example5-2} (right).}
\end{figure}

\begin{figure}[h]
  \centering
  \includegraphics[width=0.495\textwidth]{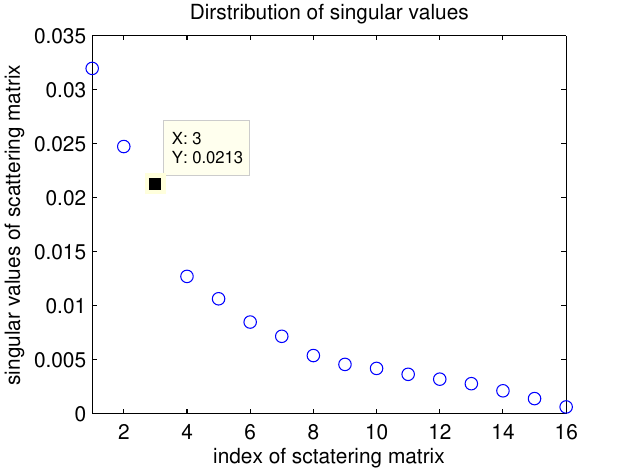}
  \includegraphics[width=0.495\textwidth]{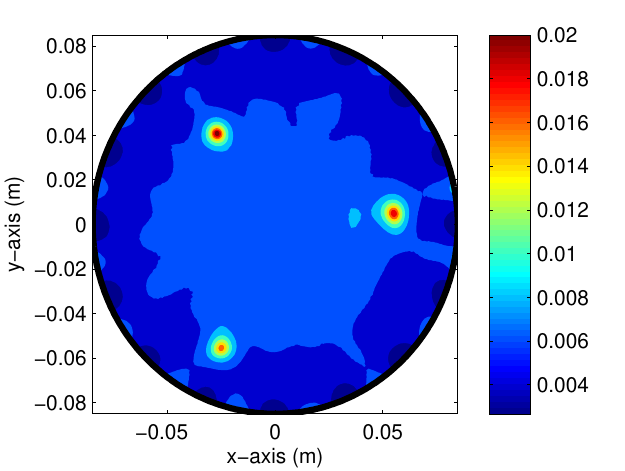}
  \caption{\label{Real1}(Example \ref{Example5-1}) Distribution of singular values of scattering matrix (left) and map of $\mathfrak{F}_{\dm}(\mr)$ (right).}
\end{figure}

\begin{ex}[Imaging of three screw drivers: small anomalies with same size]\label{Example5-1}
Figure \ref{Real1} shows the distribution of singular values of $\mathbb{D}$ and map of $\mathfrak{F}_{\dm}(\mr)$ for imaging cross-section of three screw drivers with same radii. Similar to the result in Figure \ref{CST2}, we can observe that total number of nonzero singular values is equal to the total number of small anomalies and their shapes were imaged satisfactorily.
\end{ex}

\begin{figure}[h]
  \centering
  \includegraphics[width=0.495\textwidth]{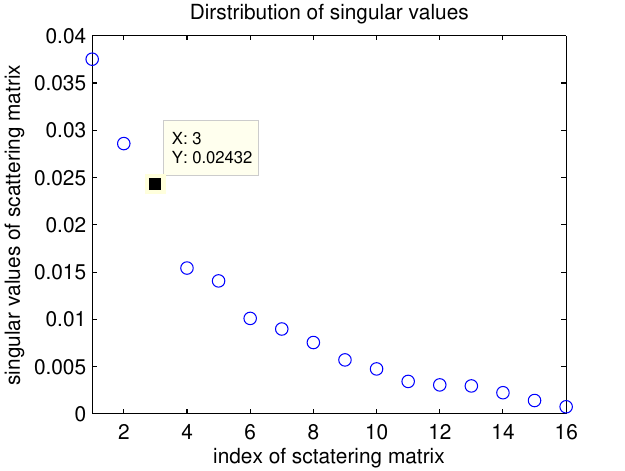}
  \includegraphics[width=0.495\textwidth]{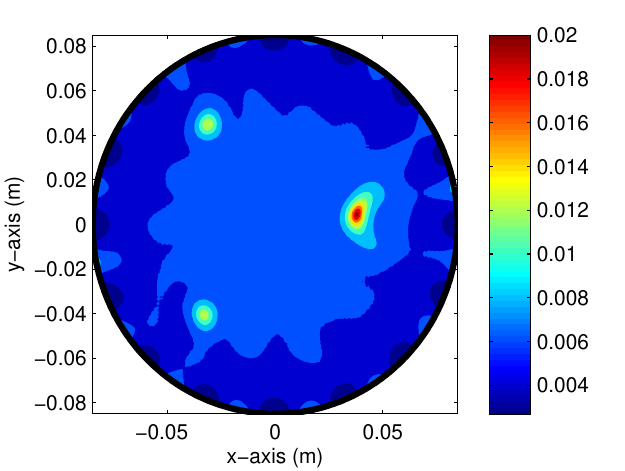}\\
  \includegraphics[width=0.495\textwidth]{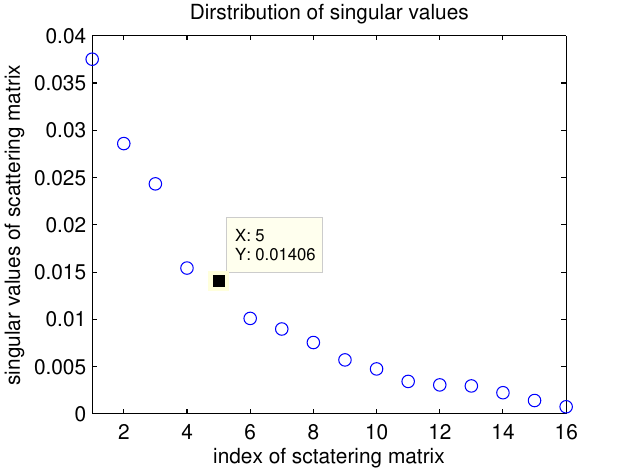}
  \includegraphics[width=0.495\textwidth]{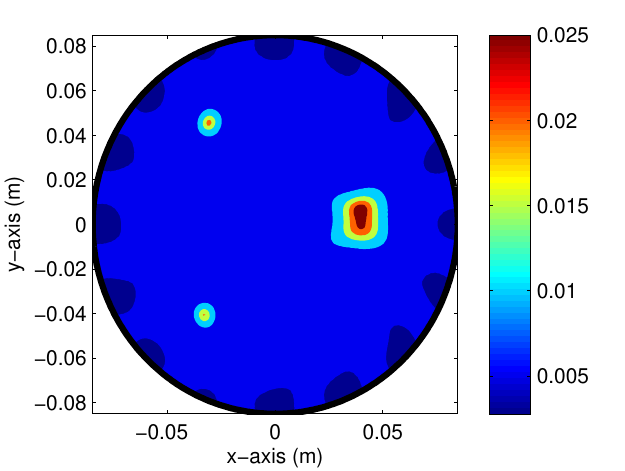}
  \caption{\label{Real2} (Example \ref{Example5-2}) Distribution of singular values of scattering matrix (left column) and map of $\mathfrak{F}_{\dm}(\mr)$ (right column).}
\end{figure}

\begin{ex}[Imaging of two screw drivers and hand hammer: small anomalies with different sizes]\label{Example5-2}
Now, we consider the imaging results of small anomalies with different sizes. For this, two screw drivers and one hand hammer were selected for cross-section imaging. Figure \ref{Real2} shows the distribution of $\mathbb{D}$ and maps of $\mathfrak{F}_{\dm}(\mr)$ with $\M=3$ and $\M=5$. Through the map of $\mathfrak{F}_{\dm}(\mr)$, it is possible to recognize that one anomaly (hand hammer) is larger than the others (drivers). However, unlike imaging of small anomalies with same size, a careful threshold to discriminate nonzero singular values is necessary for obtaining good result.
\end{ex}

\begin{figure}[h]
  \centering
  \includegraphics[width=0.495\textwidth]{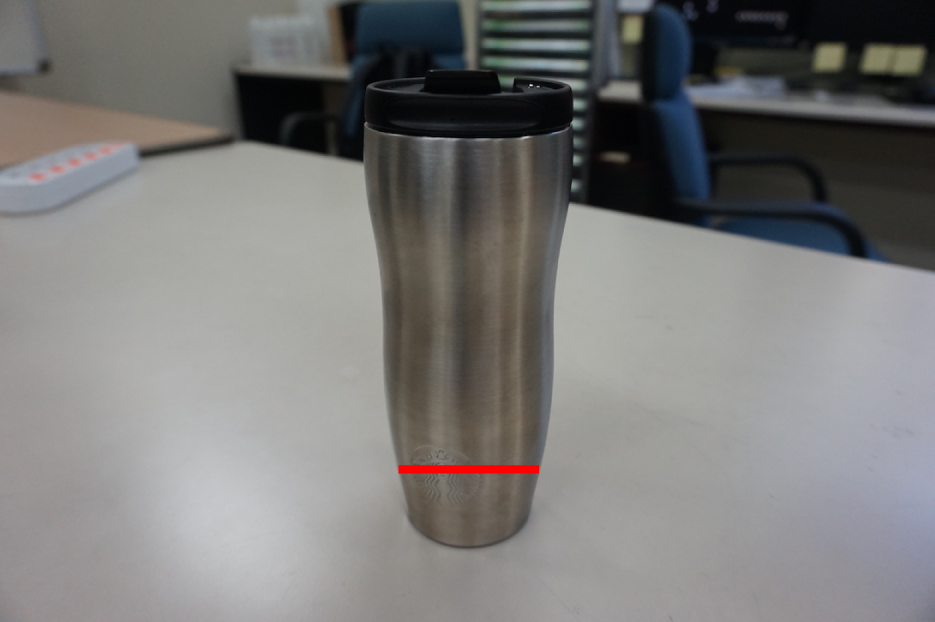}
  \includegraphics[width=0.495\textwidth]{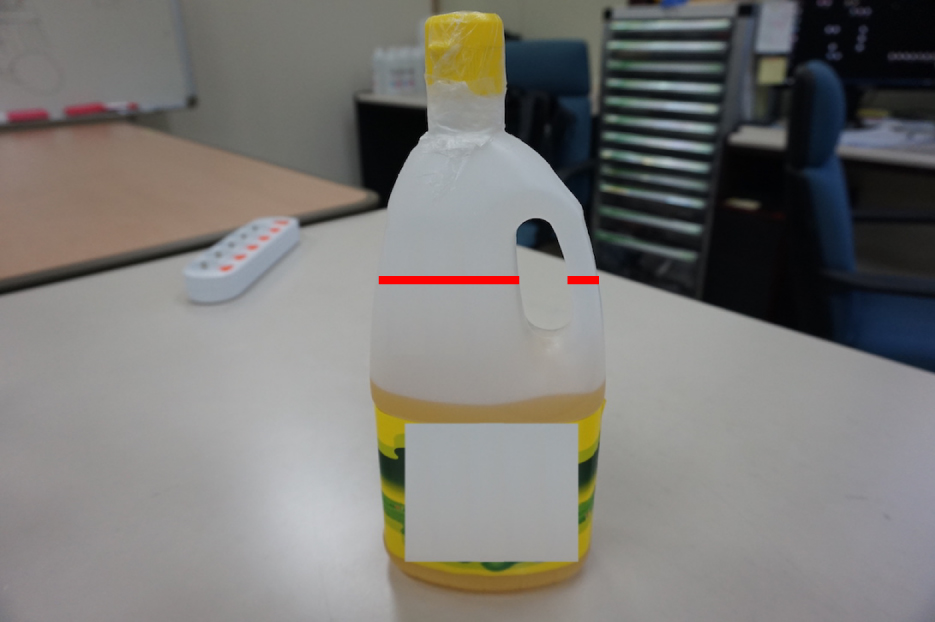}
  \caption{\label{Anomalies}Illustration of the anomalies for Examples \ref{Example5-3} (left) and \ref{Example5-4} (right). Red-colored lines describe the cross-sections of the simulation results in Figures \ref{Real3} and \ref{Real4}.}
\end{figure}

\begin{ex}[Imaging of a tumbler: extended anomaly]\label{Example5-3}
We next consider the imaging of an extended anomaly. For this, we select a large tumbler as an extended anomaly, refer to Figure \ref{Anomalies}. To obtain the map of $\mathfrak{F}_{\dm}(\mr)$ with $\M=6$ and $\M=9$ singular values were selected to define the projection operator; however, in general, we have observed that it is very hard to discriminate the nonzero singular values. Since the size of the tumbler is large, the Born approximation cannot be applied so that the it is impossible to identify complete shape of tumbler. Although it is impossible to image the complete shape, its boundary was successfully imaged. This is similar result as in \cite{HSZ1} about the imaging of an extended target with Dirichlet boundary condition. Hence, this result can be regarded as a good initial guess of Newton-type schemes or level-set strategies.
\end{ex}

\begin{figure}[h]
  \centering
  \includegraphics[width=0.495\textwidth]{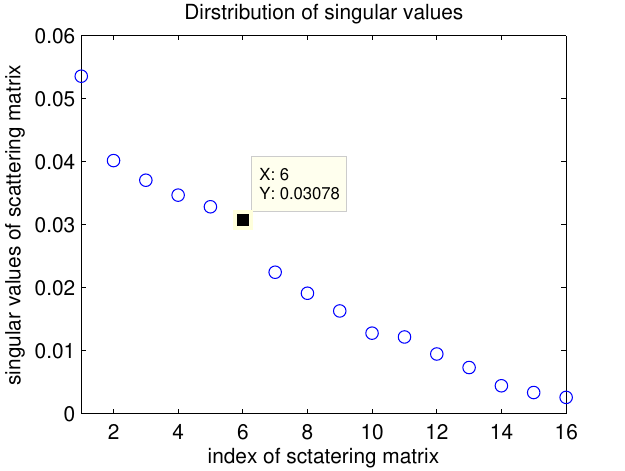}
  \includegraphics[width=0.495\textwidth]{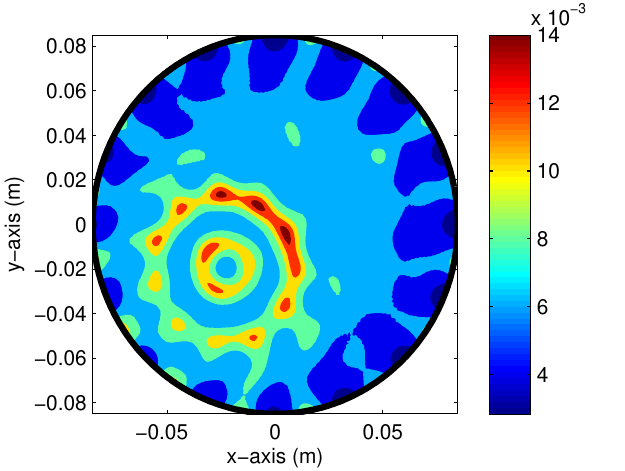}\\
  \includegraphics[width=0.495\textwidth]{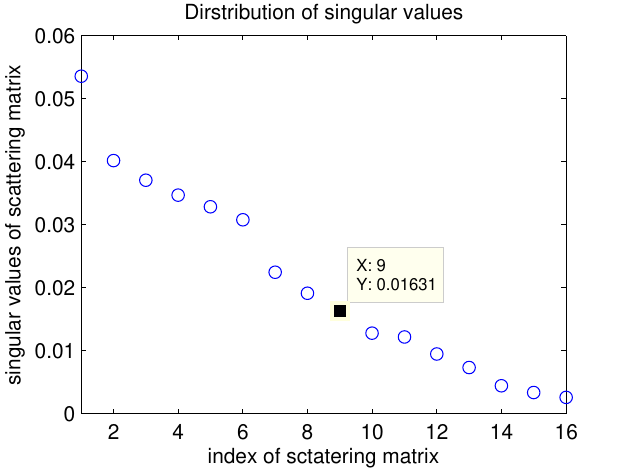}
  \includegraphics[width=0.495\textwidth]{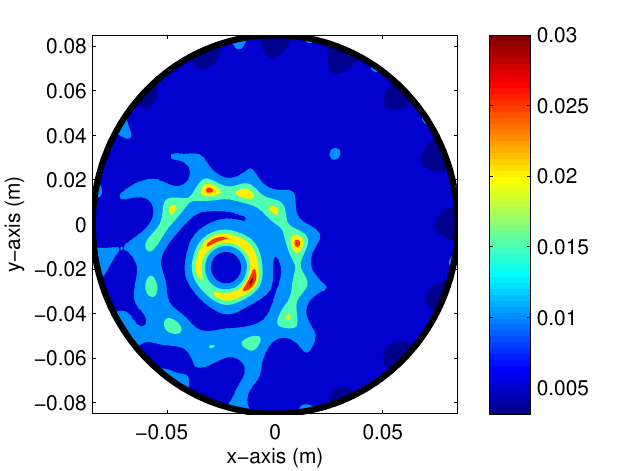}
  \caption{\label{Real3}(Example \ref{Example5-3}) Distribution of singular values of scattering matrix (left column) and map of $\mathfrak{F}_{\dm}(\mr)$ (right column).}
\end{figure}

\begin{figure}[h]
  \centering
  \includegraphics[width=0.495\textwidth]{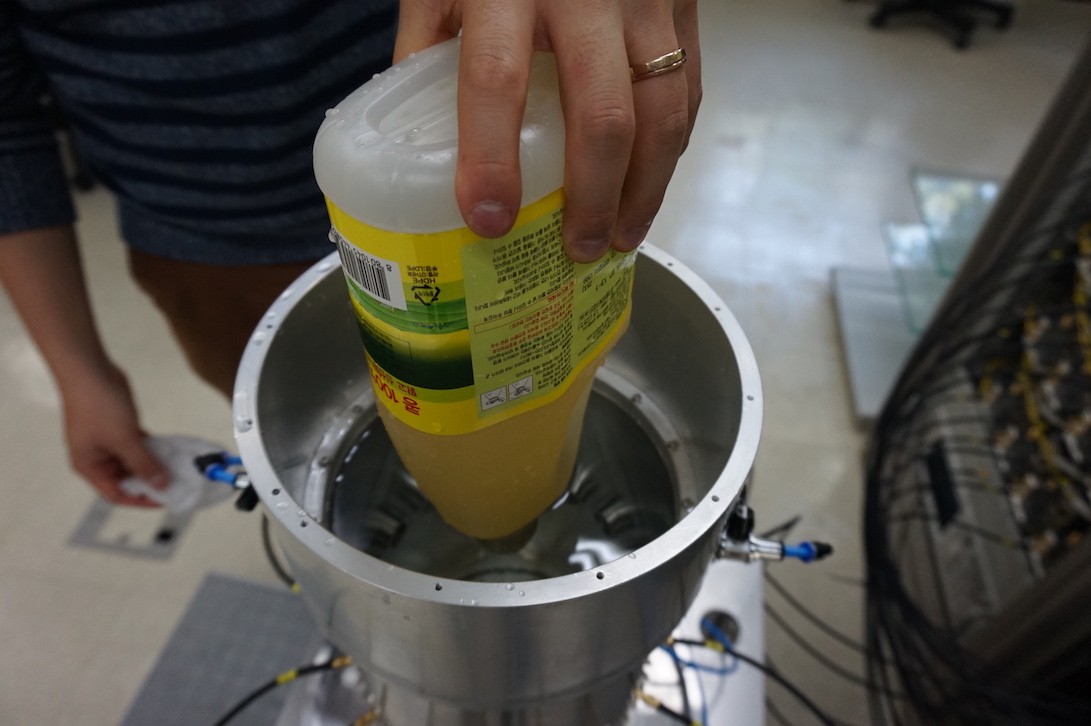}
  \includegraphics[width=0.495\textwidth]{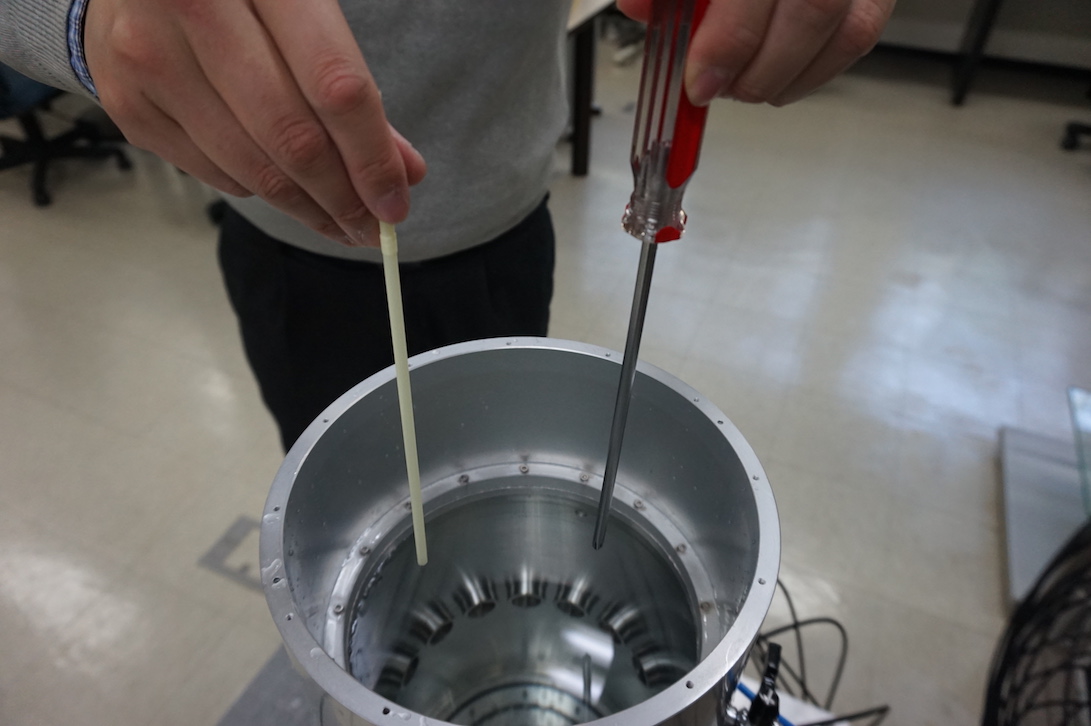}
  \caption{\label{Configuration-Real2}Real-data experiments for Examples \ref{Example5-4} (left) and \ref{Example5-5} (right).}
\end{figure}

\begin{ex}[Imaging of complex shaped anomaly: limitation of MUSIC]\label{Example5-4}
For the next example, we consider the imaging of complex shaped anomaly. For this, a plastic bottle with handgrip filled by soybean oil is chosen, refer to Figure \ref{Anomalies}. Notice that since the shape of anomaly is complex, discrimination of nonzero singular values is very hard and the only the outline shape of plastic bottle can be recognized via the map of $\mathfrak{F}_{\dm}(\mr)$ with $\M=7$. If one choose $\M=12$ for defining projection operator, the result is coarse. This result shows a limitation of application of MUSIC for a complex shaped anomaly.
\end{ex}

\begin{figure}[h]
  \centering
  \includegraphics[width=0.495\textwidth]{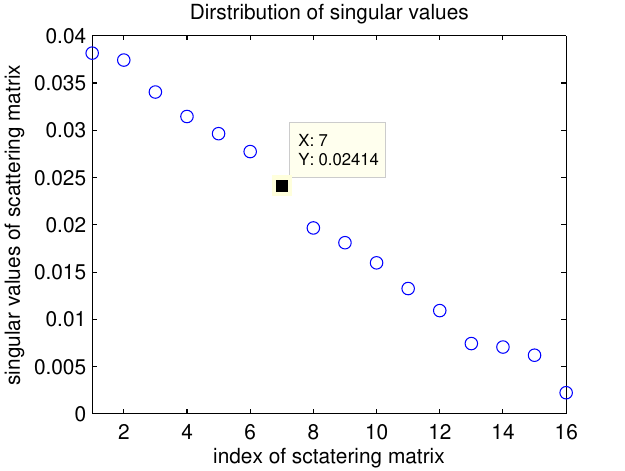}
  \includegraphics[width=0.495\textwidth]{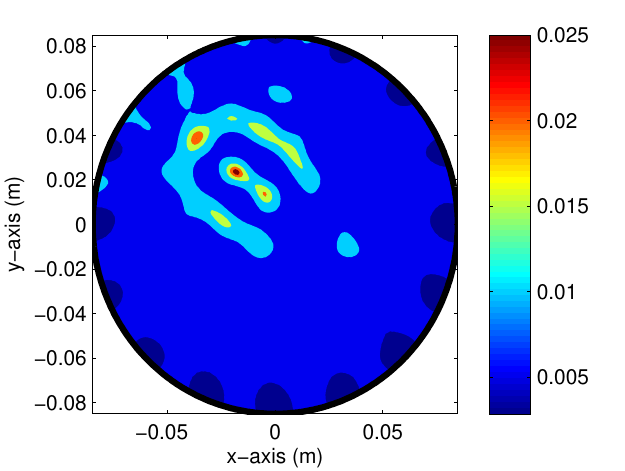}\\
  \includegraphics[width=0.495\textwidth]{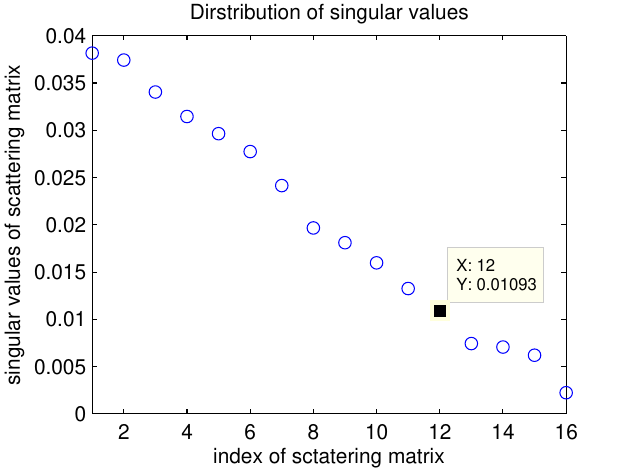}
  \includegraphics[width=0.495\textwidth]{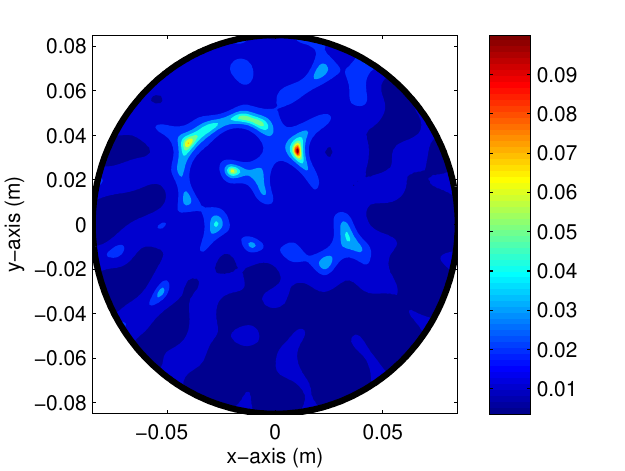}
  \caption{\label{Real4}(Example \ref{Example5-4}) Distribution of singular values of scattering matrix (left column) and map of $\mathfrak{F}_{\dm}(\mr)$ (right column).}
\end{figure}

\begin{ex}[Imaging of screw driver and plastic straw: limitation of MUSIC]\label{Example5-5}
For the final example, we exhibit a result which shows another limitation of MUSIC in real-world microwave imaging. For this, one screw driver and one plastic straw were selected. Notice that the values of permittivity of the plastic straw and screw driver are extremely small and large, respectively. This means that existence of plastic straw does not affect to the measurement data so that only one singular value of $\mathbb{D}$ that is significantly larger than the others appeared. Correspondingly, it is possible to recognize the outline shape of screw driver in the map of $\mathfrak{F}_{\dm}(\mr)$, but it is impossible to recognize the existence of plastic straw.
\end{ex}

\begin{figure}[h]
  \centering
  \includegraphics[width=0.495\textwidth]{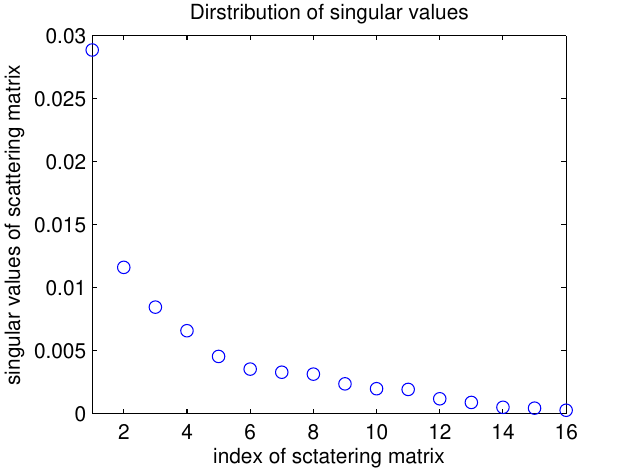}
  \includegraphics[width=0.495\textwidth]{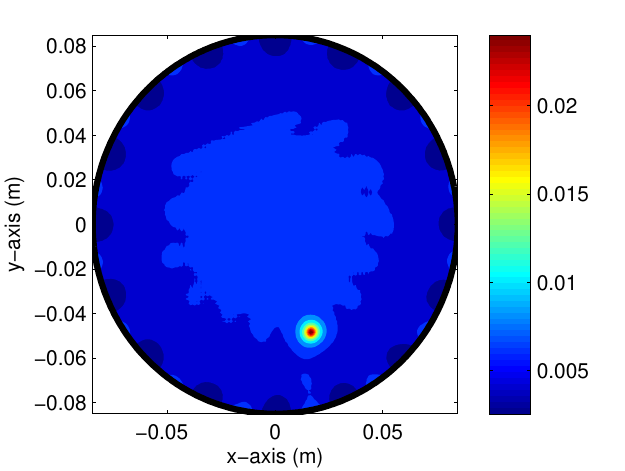}
  \caption{\label{Real5}Distribution of singular values of scattering matrix (left) and map of $\mathfrak{F}_{\dm}(\mr)$ (right).}
\end{figure}

\section{Concluding Remarks}\label{sec:6}
In this research, we applied MUSIC algorithm for a non-iterative microwave imaging of unknown anomalies when the diagonal elements of scattering matrix are unknown. To show the feasibility of the MUSIC in real-world application, we explored mathematical structure of imaging function by finding a relationship with infinite series of Bessel function of integer order. Based on explored structures, we confirmed that the shapes of small anomalies can be imaged uniquely and successfully via the MUSIC. Hence, it can be regarded as a good microwave imaging technique.

Based on several works \cite{AIL1,IGLP,AJP,JKP}, it has been shown that MUSIC is effective in a limited-aperture inverse scattering problem. Application of MUSIC in real-world limited-aperture microwave imaging will be the forthcoming work. Furthermore, as we seen, there still exists some limitations on the real-world applications. Improvement of MUSIC for obtaining better results in real-world microwave imaging will be an interesting research subject.

\section*{Acknowledgement}
The author would like to acknowledge Raffaele Solimene for his valuable comments about the MUSIC algorithm and the normalization of the test vector as well as Kwang-Jae Lee and Seong-Ho Son for helping in generating scattering parameter data from CST STUDIO SUITE and microwave machine. The author also wish to thank two anonymous referees for their comments that help to increase the quality of the paper. Part of this work was done while the author was visiting Electronics and Telecommunications Research Institute (ETRI) and Laboratoire G\'enie \'electrique et \'electronique de Paris (GeePs), Universit\'e Paris-Sud. This research was supported by the National Research Foundation of Korea (NRF) grant funded by the National Research Foundation of Korea (NRF) grant funded by the Korea government (MSIT) (No. NRF-2020R1A2C1A01005221).

\appendix

\section{Proof of the Theorem \ref{StructureImaging}}\label{sec:A}
Suppose that $\mr$ is close to $\ma_n$. Then, since
  \begin{align*}
  \lim_{\mr\to\ma_n}|H_0^{(1)}(k|\mr-\ma_n|)|&=\lim_{\mr\to\ma_n}|J_0(k|\mr-\ma_n|)+iY_0(k|\mr-\ma_n|)|\\
  &\approx\lim_{\mr\to\ma_n}\left|\frac{1}{\Gamma(1)}+i\frac{2}{\pi}\bigg(\ln(k|\mr-\ma_n|)-\ln2+\gamma\bigg)\right|\longrightarrow+\infty,
  \end{align*}
  we can observe that $|\mathbb{P}_{\noise}(\mf(\mr))|\longrightarrow+\infty$ so that $\mathfrak{F}_{\dm}(\mr)\longrightarrow+0$ as $\mr\to\ma_n$ for $n=1,2,\cdots,\N$. Here, $Y_0$ denotes the Neumann function of order zero, $\Gamma$ is the Gamma function, and $\gamma=0.5772156649\ldots$ is the Euler-Mascheroni constant (see \cite{CK} for instance).

  Now, assume that $\mr$ is not close to $\ma_n$ such that $k|\mr-\ma_n|\gg1/4$. Then, based on the asymptotic form of Hankel function {(see \cite[Theorem 2.5]{CK} for instance)},
  \begin{equation}\label{Hankel}
  -\frac{i}{4}H_0^{(1)}(k|\mr-\ma_n|)\approx\frac{1-i}{16\sqrt{k\pi}}\frac{\e^{ik|\ma_n|}}{|\ma_n|^{1/2}}\e^{-ik\frac{\ma_n}{|\ma_n|}\cdot\mr}=\frac{1-i}{16\sqrt{k\pi}}\frac{\e^{ik|\ma_n|}}{|\ma_n|^{1/2}}\e^{-ik\vt_n\cdot\mr},
  \end{equation}
  $\mW(\mr)$ of \eqref{SpanVectors} can be written by
  $\mW(\mr)\approx\Big[\e^{-ik\vt_1\cdot\mr},\e^{-ik\vt_2\cdot\mr},\cdots,\e^{-ik\vt_\N\cdot\mr}\Big]^{\mathtt{T}}$. Therefore, $\mf(\mr)$ can be represented as
  \[\mf(\mr)=\frac{\mW(\mr)}{|\mW(\mr)|}\approx\frac{1}{\sqrt{\N}}\bigg[\e^{-ik\vt_1\cdot\mr},\e^{-ik\vt_2\cdot\mr},\cdots,\e^{-ik\vt_\N\cdot\mr}\bigg]^{\mathtt{T}}.\]
Since there exists only one anomaly and no multiple scattering effect, total number of nonzero singular values is equal to one (see \cite{P-SUB11} for instance). Thus, the projection operator \eqref{Projection} is of the form
\[\mathbb{P}_{\noise}=\mathbb{I}_\N-\mU_1\mU_1^*.\]
Based on \eqref{Sparameter-Representation}, \eqref{Hankel}, and orthonormal property of singular vectors, $\mU_1\mU_1^*$ can be represented as
  \[\mU_1\mU_1^*\approx\frac{1}{\N-1}\left[
               \begin{array}{cccc}
                  0 & \e^{-ik(\vt_1-\vt_2)\cdot\mr_\star} & \cdots  & \e^{-ik(\vt_1-\vt_\N)\cdot\mr_\star} \\
                  \e^{-ik(\vt_2-\vt_1)\cdot\mr_\star} & 0 & \cdots & \e^{-ik(\vt_2-\vt_\N)\cdot\mr_\star} \\
                  \vdots & \vdots & \ddots & \vdots \\
                 \e^{-ik(\vt_\N-\vt_1)\cdot\mr_\star} & \e^{-ik(\vt_\N-\vt_2)\cdot\mr_\star} & \cdots & 0
               \end{array}
             \right].\]
With this, by defining $\mathcal{N}=\set{1,2,\cdots,\N}$, we can evaluate
\begin{align*}
(\mU_1\mU_1^*)\mf(\mr)
&\approx\frac{1}{\sqrt{\N}(\N-1)}\left[
               \begin{array}{cccc}
                  0 & \e^{-ik(\vt_1-\vt_2)\cdot\mr_\star} & \cdots  & \e^{-ik(\vt_1-\vt_\N)\cdot\mr_\star} \\
                  \e^{-ik(\vt_2-\vt_1)\cdot\mr_\star} & 0 & \cdots & \e^{-ik(\vt_2-\vt_\N)\cdot\mr_\star} \\
                  \vdots & \vdots & \ddots & \vdots \\
                 \e^{-ik(\vt_\N-\vt_1)\cdot\mr_\star} & \e^{-ik(\vt_\N-\vt_2)\cdot\mr_\star} & \cdots & 0
               \end{array}
             \right]
\left[
               \begin{array}{c}
                  \e^{-ik\vt_1\cdot\mr} \\
                  \e^{-ik\vt_2\cdot\mr} \\
                  \vdots \\
                 \e^{-ik\vt_\N\cdot\mr}
               \end{array}
             \right]\\
&=\frac{1}{\sqrt{\N}(\N-1)}\left[\begin{array}{c}
                 \medskip \displaystyle \e^{-ik\vt_1\cdot\mr_\star}\sum_{n\in\mathcal{N}\backslash\set{1}}\e^{ik\vt_n\cdot(\mr_\star-\mr)} \\
                 \medskip \displaystyle \e^{-ik\vt_2\cdot\mr_\star}\sum_{n\in\mathcal{N}\backslash\set{2}}\e^{ik\vt_n\cdot(\mr_\star-\mr)} \\
                 \medskip \vdots \\
                 \displaystyle \e^{-ik\vt_\N\cdot\mr_\star}\sum_{n\in\mathcal{N}\backslash\set{\N}}\e^{ik\vt_n\cdot(\mr_\star-\mr)}
               \end{array}
             \right].
\end{align*}
Since the following Jacobi--Anger expansion holds uniformly
\begin{equation}\label{JacobiAnger}
\e^{ix\cos\theta}=J_0(x)+\sum_{\nu\in\mathbb{Z}^*}i^\nu J_{\nu}(x)\e^{i\nu\theta},
\end{equation}
interchange of limit is available so that we can derive that for $p=1,2,\cdots,\N$,
\begin{align*}
\e^{-ik\vt_p\cdot\mr_\star}\sum_{n\in\mathcal{N}\backslash\set{p}}\e^{ik\vt_n\cdot(\mr_\star-\mr)}&=\e^{-ik\vt_p\cdot\mr_\star}\left(\sum_{n=1}^{\N}\e^{ik\vt_n\cdot(\mr_\star-\mr)}-\e^{ik\vt_p\cdot(\mr_\star-\mr)}\right)\\
&=\e^{-ik\vt_p\cdot\mr_\star}\sum_{n=1}^{\N}\e^{ik|\mr_\star-\mr|\cos(\theta_n-\phi_\star)}-\e^{-ik\vt_p\cdot\mr}\\
&=\e^{-ik\vt_p\cdot\mr_\star}\sum_{n=1}^{\N}\left(J_0(k|\mr_\star-\mr|)+\sum_{\nu\in\mathbb{Z}^*}i^\nu J_{\nu}(k|\mr_\star-\mr|)\e^{i\nu(\theta_n-\phi_\star)}\right)-\e^{-ik\vt_p\cdot\mr}\\
&:=\e^{-ik\vt_p\cdot\mr_\star}\sum_{n=1}^{\N}\bigg(J_0(k|\mr_\star-\mr|)+\mathcal{R}(\mr_\star,\ma_n)\bigg)-\e^{-ik\vt_p\cdot\mr}.
\end{align*}
Thus, $\mathbb{P}_{\noise}(\mW(\mr))$ can be written as
\begin{align*}
\mathbb{P}_{\noise}(\mW(\mr))&=\Big(\mathbb{I}_\N-\mU_1\mU_1^*\Big)\mW(\mr)\\
&\approx\frac{1}{\sqrt{\N}}\left[\begin{array}{c}
                 \medskip \displaystyle \e^{-ik\vt_1\cdot\mr}-\frac{1}{\N-1}\left\{\e^{-ik\vt_1\cdot\mr_\star}\sum_{n=1}^{\N}\bigg(J_0(k|\mr_\star-\mr|)+\mathcal{R}(\mr_\star,\ma_n)\bigg)-\e^{-ik\vt_1\cdot\mr}\right\} \\
                 \medskip \displaystyle \e^{-ik\vt_2\cdot\mr}-\frac{1}{\N-1}\left\{\e^{-ik\vt_2\cdot\mr_\star}\sum_{n=1}^{\N}\bigg(J_0(k|\mr_\star-\mr|)+\mathcal{R}(\mr_\star,\ma_n)\bigg)-\e^{-ik\vt_2\cdot\mr}\right\} \\
                 \medskip \vdots \\
                 \displaystyle \e^{-ik\vt_\N\cdot\mr}-\frac{1}{\N-1}\left\{\e^{-ik\vt_{\N}\cdot\mr_\star}\sum_{n=1}^{\N}\bigg(J_0(k|\mr_\star-\mr|)+\mathcal{R}(\mr_\star,\ma_n)\bigg)-\e^{-ik\vt_{\N}\cdot\mr}\right\}
               \end{array}
             \right]\\
&=\frac{\sqrt{\N}}{(N-1)}\left[\begin{array}{c}
                 \medskip \displaystyle \e^{-ik\vt_1\cdot\mr}-\frac{1}{\N}\e^{-ik\vt_1\cdot\mr_\star}\sum_{n=1}^{\N}\bigg(J_0(k|\mr_\star-\mr|)+\mathcal{R}(\mr_\star,\ma_n)\bigg) \\
                 \medskip \displaystyle \e^{-ik\vt_2\cdot\mr}-\frac{1}{\N}\e^{-ik\vt_2\cdot\mr_\star}\sum_{n=1}^{\N}\bigg(J_0(k|\mr_\star-\mr|)+\mathcal{R}(\mr_\star,\ma_n)\bigg) \\
                 \medskip \vdots \\
                 \displaystyle \e^{-ik\vt_\N\cdot\mr}-\frac{1}{\N}\e^{-ik\vt_{\N}\cdot\mr_\star}\sum_{n=1}^{\N}\bigg(J_0(k|\mr_\star-\mr|)+\mathcal{R}(\mr_\star,\ma_n)\bigg)
               \end{array}
             \right]
\end{align*}
and correspondingly, we arrive
\[|\mathbb{P}_{\noise}(\mW(\mr))|=\abs{\mathbb{P}_{\noise}(\mW(\mr))\cdot\overline{\mathbb{P}_{\noise}(\mW(\mr))}}^{1/2}
\approx\frac{\sqrt{\N}}{\N-1}\abs{\N-\frac{1}{\N}\sum_{m=1}^{\N}(\Phi_1+\overline{\Phi}_1)+\frac{1}{\N^2}\sum_{m=1}^{\N}\Phi_2\overline{\Phi}_2}^{1/2},\]
where
\[\Phi_1=\e^{ik\vt_m\cdot(\mr_\star-\mr)}\sum_{n=1}^{\N}\bigg(J_0(k|\mr_\star-\mr|)+\overline{\mathcal{R}(\mr_\star,\ma_n)}\bigg)\quad\text{and}\quad\Phi_2=\e^{-ik\vt_m\cdot\mr_\star}\sum_{n=1}^{\N}\bigg(J_0(k|\mr_\star-\mr|)+\mathcal{R}(\mr_\star,\ma_n)\bigg).\]

First, applying \eqref{JacobiAnger} again, we can see that
\[\e^{ik\vt_m\cdot(\mr_\star-\mr)}=J_0(k|\mr_\star-\mr|)+\sum_{\nu\in\mathbb{Z}^*}i^\nu J_{\nu}(k|\mr_\star-\mr|)\e^{i\nu(\theta_m-\phi_\star)}=J_0(k|\mr_\star-\mr|)+\mathcal{R}(\mr_\star,\ma_m).\]
With this, performing an elementary calculus yields
\begin{align*}
\sum_{m=1}^{\N}\Phi_1&=\left(\N J_0(k|\mr_\star-\mr|)+\sum_{n=1}^{\N}\mathcal{R}(\mr_\star,\ma_m)\right)\left(\N J_0(k|\mr_\star-\mr|)+\sum_{n=1}^{\N}\overline{\mathcal{R}(\mr_\star,\ma_n)}\right)\\
&=\N^2J_0(k|\mr_\star-\mr|)+\N J_0(k|\mr_\star-\mr|)\sum_{n=1}^{\N}\bigg(\mathcal{R}(\mr_\star,\ma_n)+\overline{\mathcal{R}(\mr_\star,\ma_n)}\bigg)+\left(\sum_{n=1}^{\N}\mathcal{R}(\mr_\star,\ma_n)\right)\left(\sum_{n=1}^{\N}\overline{\mathcal{R}(\mr_\star,\ma_n)}\right).
\end{align*}
Hence, we can obtain
\begin{multline}\label{term1}
\frac{1}{\N}\sum_{m=1}^{\N}(\Phi_1+\overline{\Phi}_1)=2\N J_0(k|\mr_\star-\mr|)^2+4J_0(k|\mr_\star-\mr|)\text{Re}\sum_{n=1}^{\N}\mathcal{R}(\mr_\star,\ma_n)\\
+\frac{1}{\N}\text{Re}\left(\sum_{n=1}^{\N}\mathcal{R}(\mr_\star,\ma_n)\right)\left(\sum_{n=1}^{\N}\overline{\mathcal{R}(\mr_\star,\ma_n)}\right).
\end{multline}

Next, by performing an elementary calculus,
\begin{align*}
\Phi_2\overline{\Phi}_2&=\left(\N J_0(k|\mr_\star-\mr|)+\sum_{n=1}^{\N}\mathcal{R}(\mr_\star,\ma_n)\right)\left(\N J_0(k|\mr_\star-\mr|)+\sum_{n=1}^{\N}\overline{\mathcal{R}(\mr_\star,\ma_n)}\right)\\
&=\N^2J_0(k|\mr_\star-\mr|)+\N J_0(k|\mr_\star-\mr|)\sum_{n=1}^{\N}\bigg(\mathcal{R}(\mr_\star,\ma_n)+\overline{\mathcal{R}(\mr_\star,\ma_n)}\bigg)+\left(\sum_{n=1}^{\N}\mathcal{R}(\mr_\star,\ma_n)\right)\left(\sum_{n=1}^{\N}\overline{\mathcal{R}(\mr_\star,\ma_n)}\right).
\end{align*}
Hence, we can derive
\begin{equation}\label{term2}
\frac{1}{\N^2}\sum_{m=1}^{\N}\Phi_2\overline{\Phi}_2=\N J_0(k|\mr_\star-\mr|)+2J_0(k|\mr_\star-\mr|)\text{Re}\sum_{n=1}^{\N}\mathcal{R}(\mr_\star,\ma_n)+\frac{1}{\N}\text{Re}\left(\sum_{n=1}^{\N}\mathcal{R}(\mr_\star,\ma_n)\right)\left(\sum_{n=1}^{\N}\overline{\mathcal{R}(\mr_\star,\ma_n)}\right).
\end{equation}

Finally, combining \eqref{term1} and \eqref{term2}, we can see that
\begin{multline*}
|\mathbb{P}_{\noise}(\mW(\mr))|=\frac{\sqrt{\N}}{\N-1}\left|\N-\N J_0(k|\mr_\star-\mr|)^2+2J_0(k|\mr_\star-\mr|)\text{Re}\sum_{n=1}^{\N}\mathcal{R}(\mr_\star,\ma_n)\right.\\
+\left.\frac{1}{\N}\text{Re}\left(\sum_{n=1}^{\N}\mathcal{R}(\mr_\star,\ma_n)\right)\left(\sum_{n=1}^{\N}\overline{\mathcal{R}(\mr_\star,\ma_n)}\right)\right|^{1/2}
\end{multline*}
and the structure \eqref{StructureImagingFunction} can be derived correspondingly. This completes the proof.

\end{document}